\documentclass[10pt]{article}

\usepackage{amssymb}
\usepackage{latexsym}
\usepackage{epsfig}

\newcommand{\mod}[1]{{\rm\ (mod\ }#1{\rm)}}

 \newcommand{\RR}{\mathbb R}
 
 \newcommand{\ZZ}{\mathbb Z}

\newtheorem{theorem}{Theorem}
\newtheorem{lemma}{Lemma}

\newtheorem{coro}{Corollary}

 \newcommand{\be}{\begin{equation}}
 \newcommand{\ee}{\end{equation}}

\begin{document}
 \bibliographystyle{unsrt}

\title{Lattice Substitution Systems and Model Sets}

\author{{\sc Jeong-Yup Lee and Robert V. Moody}\\
Department of Mathematical Sciences, \\
University of Alberta,\\
Edmonton, T6G 2G1, Canada}

\maketitle

{\footnotesize
\flushright{
Learning is but an adjunct to ourself\\
And where we are our learning likewise is.\\
- W. Shakespeare\\
}}

\begin{abstract}\par
The paper studies ways in which the sets of a partition of a lattice in
$\RR^n$ become regular model sets. The  main theorem
gives equivalent conditions which assure
that a matrix substitution system on a lattice in $\RR^n$ gives rise to 
regular model sets (based on $p$-adic-like internal spaces), and hence
to pure point diffractive sets.
The methods developed here are used to show
that the $n-$dimensional chair tiling and the sphinx tiling are
 pure point diffractive.
\end{abstract}

\section{Introduction}
There have been two very successful approaches to building discrete mathematical
structures with long-range aperiodic order. These are the substitution methods,
notably symbolic substitutions and tiling substitutions, and the cut and project
method. In the first case the structure is typically generated by successive substitution
from a finite starting configuration. In the second it typically appears in one shot as the 
(partial) projection of a periodic structure in some ``higher'' dimensional embedding
space.

The principal focus in this paper is the relationship between matrix
substitution systems on a {\em lattice} and a naturally related cut and
project formalism. We start with a partition of a lattice $L$ in
$\RR^n$ into a finite number of point sets $\tilde{U} = (U_1, \dots,
U_m)$ and a finite set of substitution rules $\Phi$ which are affine
inflations and under which $\tilde{U}$ is invariant. The main theorem
(Theorem \ref{mainTheorem}) provides conditions on $\Phi$ which are equivalent to $U_1,
\dots , U_m$ being regular model sets (i.e. cut and project sets). 
One of the characterizations (modular coincidence) affords a simple computational
approach to testing for model sets. In a later section
we go beyond the context of substitution systems  and provide an alternative characterization
(Theorem \ref{mainCITheorem}) of model sets. We use both types of characterization 
in showing that the sphinx and $n$-dimensional chair tilings are based on model sets. 

The connection between substitution systems and cut and project sets is
nothing new, e.g. the Fibonacci chain is often described in terms of a
cut by a strip through $\ZZ^2$, and the klotz construction of Kramer et
al. \cite{michael} is a sophisticated elaboration of the same idea.
Nonetheless, substitution systems and cut and project sets are not
different formulations of the same thing, and the relationship between
them remains inadequately understood.

In the early study of aperiodic order, the cut and project formalism was always based on
projection into $\RR^n$ from a lattice in some higher space $\RR^n \times \RR^p$, the
projection being controlled by a compact set $W \subset \RR^p$. However, it was already
implicit in the much earlier work of Y. Meyer \cite{YM1} that $\RR^p$ can be replaced
by any locally compact abelian group $H$ and $W \subset H$ by any compact set with non-empty
interior, and the projection method still produces discrete aperiodic sets with diffractive
properties (hence long-range order). Meyer's terminology for such sets 
was ``model sets'' and we use it
here in deference to its priority and to emphasize the greater generality of the internal space
$H$. Model sets have been studied in detail in 
\cite{Lagarias1,YM2,RVM1, RVM2, martin2, martin3}. 
The relevance of more
general internal spaces to tiling theory and symbolic substitutions was made explicit in
\cite{padics} where $p$-adic and mixed $p$-adic and real spaces naturally appear. 

One of the important features of making the connection to model sets is
that once it is established, pure point diffractivity is assured (see
Theorem \ref{martinThm} for a precise statement of this). This type of information is
generally quite difficult to obtain. For example, our results shown here prove
that the $n$-dimensional chair tiling and the $2$-dimensional sphinx
tiling are pure point diffractive. The former has been established for
$n=2$ previously \cite{padics,boris}. The latter is claimed in
\cite{boris} as being provable 
by a geometric form of ``coincidence'' established there (see below for more
on the concept of coincidence).

The $p$-adic type internal spaces occur when the aperiodic set in question is based
on the points of a lattice and its sublattices in $\RR^n$. An important class of 
examples of this type arises from the {\em equal length} symbolic substitution systems.
Suppose that $A =\{a_1, \dots ,a_m\}$ is a finite alphabet with associated monoid
of words $A^*$, and we are given a primitive substitution $\sigma:A \longrightarrow A^*$
for which the length $l$ of each of the words $\sigma(a_i)$ is the same. This substitution
leads to a tiling of $\RR$ of tiles of equal length, say equal to $1$.
 Matching the coordinate of the left end of
each tile with its tile type $a_i$, we obtain a partition $U_1 \cup \dots \cup U_m$
of $\ZZ$, and $\sigma$ may be viewed as comprised of a set of affine mappings
$x \mapsto lx+v$ where $v \in \ZZ$.

A lot more is known about equal length substitutions than the arbitrary ones, a 
particularly important example of this being Dekking's criterion for diffraction
\cite{dekking}. An equal length aperiodic tiling is pure point diffractive if and only if
it admits a coincidence ($\sigma$ is said to admit a coincidence if there is
a $k$, $1 \le k \le l^n$,  for which the $k$th letter of each word 
$\sigma^n(a_i)$ for some $n$ is the same).

In this paper we prove a related result, but this time the dimension is not restricted.
Namely, there is a notion of coincidence (in fact there are two such notions) and
either of these is equivalent to the sets $U_1, \dots, U_m$ being regular model sets. 
One of the criteria for coincidence that we give is a straightforward algorithm
and thus in principle is computable.

As we have already pointed out, a consequence of our result is that coincidence implies 
the pure point
diffractivity  of $U_1, \dots, U_m$. 
We do not know yet to what extent the condition is equivalent to pure point diffractivity. 

The setting of the paper is entirely at the level of point sets, so necessarily the
strong conditions implicit in the tiling situation are replaced here by a corresponding
algebraic condition
on the matrix substitution system: the Perron-Frobenius eigenvalue of 
the substitution system should equal its inflation constant. This is in fact a compatibility
condition which is necessary for the model set connection to exist. 
This condition, not surprisingly, has occurred elsewhere in the literature (see 
for instance, \cite{LW,boris}). The important
result that gets the process off the ground is Theorem \ref{mainPFtheorem}, which is
largely due to Martin Schlottmann.

Matrix substitution systems, treated at the level of point
sets, have recently appeared in Lagarias and Wang \cite{LW} under
the name of self-replicating Delone sets. In that paper, point sets $X$
are not restricted to lattices and the principal question revolves
around the interesting question of existence of tilings of $\RR^n$  
by translations of certain prototiles for which the points of $X$ are 
the appropriate translational vectors.
 Also related to our paper is the study of sets of affine mappings in the context 
of lattice tilings (see \cite{vince} for a nice recent survey on this). In relation to our 
paper, the situation there corresponds to the $1\times 1$ matrix substitution systems and
the problems become entirely different. Since the tilings there are lattice tilings,
the whole issue of model sets and diffracion is trivial, and the issues lie more
around the complex nature of the tiles themselves. 

\section{Definitions and Notation}

Let $X$ be a nonempty set. 
For $m \in \ZZ_{+}$, an $m \times m$ {\em matrix function system} (MFS)
 on $X$ is an $m \times m$ matrix $\Phi = (\Phi_{ij})$, where each $\Phi_{ij}$ 
is a set (possibly empty) of mappings $X$ to $X$.

The corresponding matrix $S(\Phi) := (\mbox{card} (\Phi_{ij}))_{ij}$ is called 
the {\em substitution matrix} of $\Phi$.  The MFS is {\em primitive} if $S(\Phi)$ 
is primitive, i.e. there is an $l > 0$ for which $S(\Phi)^{l}$ has no zero entries.  

In this paper we deal only with MFSs which are {\em finite} in the sense 
that card$(\Phi_{ij}) < \infty$  for all $i,j$.  Of particular importance are the
Perron-Frobenius (PF) eigenvalue and the corresponding PF eigenvector (unique up to a scalar factor)
of $S(\Phi)$.
We will also have use for the {\em incidence matrix} $I(\Phi)$ of $\Phi$, 
which is defined by 
\[ (I(\Phi))_{ij} = \left\{ \begin{array}{ll}
                           1  & ~~~ \mbox{if}~~~ \mbox{card}(\Phi_{ij}) \neq 0 \, ,\\
                           0  & ~~~ \mbox{else}\,. 
                           \end{array}
                    \right. \]
   
Let $P(X)$ be the set of subsets of $X$.
Any MFS induces a mapping on $P(X)^{m}$ by
\begin{eqnarray}
{\Phi \left[ \begin{array}{c}
                     U_{1} \\ \vdots \\ U_{m}
             \end{array}
      \right]}&=&{\left[ \begin{array}{c}             
                         \bigcup_{j} \bigcup_{f \in \Phi_{1j}} f(U_{j})\\ \vdots \\
 \bigcup_{j} \bigcup_{f \in \Phi_{mj}} f(U_{j})
                         \end{array}                    
                  \right]_{\,,}}
\end{eqnarray}          
which we call the {\em substitution} determined by $\Phi$.
We sometimes write $\Phi_{ij} (U_{j})$ to mean $\bigcup_{f \in \Phi_{ij}} f(U_{j})$.
   
In the sequel, $X$ will be a lattice $L$ in $\RR^{n}$ and the mappings of $\Phi$ will 
always be affine linear mappings of the form $x \mapsto Qx + a$, 
where $Q \in \mbox{End}_{\ZZ}(L)$ is the {\em same} for all the maps.
 Such maps extend to $\RR^{n}$. 
For any
 affine mapping $f : x \mapsto Qx + b$ on $L$ we denote the translational 
part,\,$b$,\,of $f$ by $t(f)$. 
We say that $f, g \in \Phi$ are {\em congruent} mod $QL$ if $t(f) \equiv t(g)$ mod $QL$.
 This equivalence relation partitions $\Phi$ into congruence classes. 
 For $a \in L$, 
$\Phi[a] := \{\,f \in \bigcup_{i,j}\Phi_{ij}\,|\,t(f) \equiv a \; \mbox{mod}\,QL\,\}.$ 
  
We say that $\Phi$ {\em admits a coincidence} if there is an  $i,~ 1 \leq i \leq m$,  
for which $\bigcap_{j=1}^{m} \Phi_{ij} \neq \emptyset$, i.e. the same map appears 
in every set of the $i$-th row for some $i$. Furthermore, if $\Phi^M[a]$
is contained entirely in one row of the MFS $(\Phi^M)$ for some $M > 0$,
$a \in L$, then we say that $(\tilde{U}, \Phi)$ admits a {\em modular
coincidence}.

Let $\Phi, \Psi$ be $m \times m$ MFSs on $X$.
Then we can compose them :
\be 
\Psi \circ \Phi = ((\Psi \circ \Phi)_{ij})\,,
\ee
where  $(\Psi \circ \Phi)_{ij} = \bigcup_{k=1}^{m} \Psi_{ik} \circ 
\Phi_{kj}~$and$~\Psi_{ik} \circ \Phi_{kj} := 
\left\{ \begin{array}{l} \{~g \circ f~ |~ g \in \Psi_{ik}, f \in \Phi_{kj}~\} \\
 \,\emptyset ~~~~~~ \mbox{if}~ \Psi_{ik} = \emptyset~ \mbox{or}~ \Phi_{kj} = \emptyset
\,.                      \end{array}
  \right.$
Evidently, $S(\Psi \circ \Phi) \leq S(\Psi)\, S(\Phi)$ (see  (\ref{inequal}) for the definition
of the partial order).

For an $m \times m$ MFS $\Phi$,  we say that  $~\tilde{U} := [U_{1},\cdots,U_{m}]^{T} \in P(X)^{m}$ 
is a {\em fixed point} of $\Phi$  if  $\Phi\, \tilde{U} = \tilde{U}$.

\section{Substitution Systems on Lattices}

Let $L$ be a lattice in $\RR^{n}$. A mapping $Q \in \mbox{End}_{\ZZ}(L)$ is 
an {\em inflation} for $L~$ if  $\det Q \ne 0$ and 
\be \label{inflation}
\bigcap_{k=0}^{\infty} Q^{k}L = \{0\}.
\ee
Let $Q$ be an inflation. Then $q := |\det Q| = [L : QL] > 1$.
We define the $Q$-adic completion 
\be
 \overline{L} = \overline{L_{Q}} = \lim_{\leftarrow k}~ L/Q^{k}L 
\ee
of $L$. $\overline{L}$ will be supplied with the usual topology of a profinite group. 
In particular, the cosets $a + Q^{k}\overline{L},~ a \in L,~ k = 0,1,2,\cdots$\,, 
form a basis of open sets of $\overline{L}$ and each of these cosets is both open and closed. 
When we use the word {\em coset} in this paper, we mean either a coset of the form $a +
Q^{k}\overline{L}$ in $\overline{L}$ or $a + Q^{k}L$ in $L$, according to the context. An important
observation is that any two cosets in $\overline{L}$ are either disjoint or one is contained in the other.
The same applies to cosets of $L$.
   
We let $\mu$ denote Haar measure on $\overline{L}$, normalized so 
that $\mu(\overline{L}) = 1$. Thus for cosets,
\be
 \mu(a + Q^{k}\overline{L}) = \frac{1}{|\det Q|^{k}} = \frac{1}{q^{k}}\,.
\ee
We also have need of the metric $d$ on $\overline{L}$ defined via the standard norm:
\be \|x\| := \frac{1}{q^{k}}~~~~ \mbox{if}~~ x \in Q^{k}\overline{L}~
 \backslash~ Q^{k+1}\overline{L}, ~~~\|0\| = 0\,.
\ee

From $\bigcap_{k=0}^{\infty}Q^{k}L = \{0\}$, we conclude that the mapping $x
\mapsto \{ x~ \mbox{mod}~ Q^{k}L\}_k$ embeds $L$ in $\overline{L}$. 
We identify $L$ with its image in $\overline{L}$. Note that $\overline{L}$ is the closure of $L$,
whence the notation.

{\em An affine lattice substitution system on $L$ with inflation $Q$} is a 
pair $(\tilde{U},\Phi)$ consisting of disjoint
subsets $\{U_{i}\}_{i=1}^{m}$ of $L$ and an $m \times m$ MFS $\Phi~$ 
on $L$ for 
which $\tilde{U} = [U_{1},\cdots,U_{m}]^{T}$ is a fixed point of
$\Phi$, i.e.
\be 
U_{i} = \bigcup_{j=1}^{m} \bigcup_{f \in \Phi_{ij}} f(U_{j}), ~~i = 1,\cdots,m\,, \label{8}
\ee
where the maps of $\Phi$ are affine mappings of the form $x \mapsto Qx + a,~ a \in L~$, 
and in which the unions in (\ref{8}) are {\em disjoint}.\footnote{
In the case that one has unions (\ref{8}) which are not disjoint there arises
the natural question of the mulitplicities of points, or more generally
densities of points. For more on this see \cite{BM, LW}.} In this paper
all our matrix substitution systems are composed of affine mappings on a
lattice
and we often drop the words `affine lattice', speaking simply of 
substitution systems.
   
We say that the substitution system $(\tilde{U},\Phi)$ is {\em primitive} if $\Phi$ is 
primitive. A second substitution system $(\tilde{U^{'}},\Psi)$ is called 
{\em equivalent} 
to $(\tilde{U},\Phi)$ if $\tilde{U^{'}} = \tilde{U},\,$ $\Psi$ and $\Phi$ 
have the same inflation, 
and $S(\Psi),S(\Phi)$ have the same PF-eigenvalue and right PF-eigenvector 
(up to scalar factor).
  
Let $(\tilde{U},\Phi)$ be a substitution system on $L$. Identifying $L$ as a dense 
subgroup of $\overline{L}$, 
we have a unique extension of $\Phi$ to a MFS on $\overline{L}$ in the obvious way. 
Thus if $f \in \Phi_{ij}$ and $f : x \mapsto Qx + a$,  then this formula defines a mapping 
on $\overline{L}$, to which we give the same name. Note that $f$ is a 
{\em contraction} on $\overline{L}$, since $\|Qx\| = \frac{1}{q}\|x\|\,$ for all 
$x \in \overline{L}$. Thus $\Phi$ determines a multi-component iterated function 
system on $\overline{L}$. 
Furthermore defining the compact subsets
\be
 W_{i} := \overline{U_{i}}, ~~ i = 1,\cdots,m\,, 
\ee
and using (\ref{8}) and the continuity of the mapping, we have
\be
 W_{i} = \bigcup_{j=1}^{m} \bigcup_{f \in \Phi_{ij}} f(W_{j}), ~~i = 1,\cdots,m \,, \label{10}
\ee
which shows that $\tilde{W} = [W_{1},\cdots,W_{m}]^{T}$ is the unique {\em attractor} of $\Phi$
(see \cite{BM, Hutchinson}).
  
We call $(\tilde{W},\Phi)$ the {\em associated Q-adic system}. We cannot expect in 
general that the decomposition in (\ref{10}) will be disjoint, so we do {\em not} 
call $(\tilde{W},\Phi)$ a substitution system.\\

For $X,Y \in \RR^{n}$, we write 
\[
\begin{array}{c}
 X \leq Y ~~~\mbox{if}~~ X_{i} \leq Y_{i} ~~~\mbox{for all}~~ 1 \leq i \leq n \\
 X < Y ~~~\mbox{if}~~ X_{i} < Y_{i}~ ~~\mbox{for all}~~ 1 \leq i \leq n\,.
\end{array}
\]
Similarly,   for $A,B \in M_n(\RR)$
\be \label{inequal}
\begin{array}{c}
 A \leq B ~~~\mbox{if}~~  A_{ij} \leq B_{ij} ~~~\mbox{for all}~~ 1 \leq i,j \leq n \\
 A < B ~~~\mbox{if}~~ A_{ij} < B_{ij} ~~~\mbox{for all}~~ 1 \leq i,j \leq n \,.
\end{array}
\ee

We begin by recalling a couple of results from the Perron-Frobenius theory. 

\begin{lemma} \label{PF1}
   Let $A$ be a non-negative primitive matrix with PF-eigenvalue $\lambda$.
If $0 \leq \lambda X \leq AX$, then $AX = \lambda X$.
\end{lemma}
{\sc{proof}} :   We can assume $X \neq 0$. 
Since $0 \leq \lambda X$ and $\lambda > 0$, $X \geq 0$. 
Let $X' > 0$ be a PF right-eigenvector of $A$.
Let $\alpha =\mbox{max} \{~\frac{X_{i}}{X'_{i}}~|~ 1 \leq i \leq m~\}$. 
Then $X \leq \alpha X'$ and $X$ is not strictly less than $\alpha X'$.
Claim $X = \alpha X'$. If $X \neq \alpha X'$, 
then $0 < A^{N}(\alpha X' - X) = \alpha \lambda^{N} X' - A^{N} X$ ~ for 
some $N$, since A is primitive. 
So $\lambda^{N} X \leq A^{N} X < \alpha \lambda^{N} X',$  i.e. $X < \alpha X'$. 
This is a contradiction. Therefore $AX = \lambda X$. \hfill $\square$
\begin{lemma} \label{PF2}
  Let $\lambda$ be the PF-eigenvalue of the non-negative primitive matrix $A$
 and $\mu$ be an eigenvalue of 
a matrix $B$ where $0 \leq B \leq A$. If $A \neq B$, then $|\mu| < \lambda$.
\end{lemma}
{\sc{proof}} :  Let $Y$ be a right eigenvector for eigenvalue $\mu$ of $B$, 
with $Y = [Y_{1},\cdots,Y_{m}]^T.$ 
Let $\overline{Y} = [|Y_{1}|,\cdots,|Y_{m}|]^{T} \neq 0$. 
Then $|\mu| \overline{Y} \leq B \overline{Y} \leq A
\overline{Y}$.  Let $\overline{X}^{T}$ be a positive left eigenvector for $A$ with PF-eigenvalue
$\lambda$.  So $|\mu| \overline{X}^{T} \overline{Y} \leq \overline{X}^{T} B \overline{Y} \leq 
\overline{X}^{T} A \overline{Y} = \lambda \overline{X}^{T} \overline{Y}$. 
This shows that $|\mu| \leq \lambda$. If $|\mu| = \lambda$, 
then $\lambda \overline{Y} \leq A \overline{Y}$. 
By Lemma \ref{PF1}, $\lambda \overline{Y} = A \overline{Y}$. Since $A$ is a primitive matrix, 
$\lambda^{m} \overline{Y} = A^{m} \overline{Y} > 0$ for some $m$. 
So $\overline{Y} > 0$. From $\lambda \overline{Y} \leq B \overline{Y} 
\leq A \overline{Y} = \lambda \overline{Y}$,  we have $A \overline{Y} = B \overline{Y}$. 
Therefore $A = B$. \hfill $\square$

\begin{lemma}
  Let $(\tilde{U},\Phi)$ be a primitive substitution system. 
Then for all $l = 1,2,\cdots, ~(\tilde{U}, \Phi^{l})$ is a primitive substitution system.
\end{lemma}
{\sc{proof}} : Let $i,j,k \in \{1,2,\cdots,m\}$. All the maps $g \in \Phi_{ik}~$
have domain $U_{k}$ and disjoint images in $U_{i}$. Moreover all the mappings
$g$ are injective. Likewise all the maps $f$
of
$\Phi_{kj}$ have domain $U_{j}$  and disjoint images in $U_{k}$. Thus all the maps $g \circ f \in
\Phi_{ik} \circ
\Phi_{kj}$  have domain $U_{j}$ and disjoint images in $U_{i}$. 
Furthermore $\Phi^{2} \tilde{U} = \Phi(\Phi \tilde{U}) = \Phi(\tilde{U}) = \tilde{U}$. 
So $(\tilde{U},\Phi^{2})$ is a substitution system. The argument extends in the same 
way to $(\tilde{U},\Phi^{l})$. The statement on primitivity is clear. \hfill $\square$

\begin{theorem} \label{mainPFtheorem}
  Let $(\tilde{U},\Phi)$ be a primitive substitution system with inflation $Q$ on $L$.
Let $(\tilde{W},\Phi)$ be the corresponding associated $Q$-adic system. 
Suppose that the PF-eigenvalue of $S(\Phi)$ is $|\det Q|$ 
and $\overline{L} = \bigcup_{i=1}^{m} W_{i}$.  Then 
\[
\begin{array}{l} \mbox{(i)}~~ S(\Phi^{r}) = (S(\Phi))^{r},\,r \geq 1;\\
                 \mbox{(ii)}~~ \mu(W_{i}) = 
\frac{1}{q^{r}} \sum_{j=1}^{m} (S(\Phi^{r}))_{ij} \mu(W_{j}),~\mbox{for all}\;\; 
i = 1,\cdots,m,~ r \geq 1;\\   
                 \mbox{(iii)}~~ \mbox{For all}\; \;i = 1,\cdots,m\,,~\stackrel{\circ}{W_{i}} 
                     \neq \emptyset~\mbox{and}~\mu (\partial W_{i}) = 0.
\end{array}
\]
\end{theorem} 
{\sc{proof}} : For every measurable set $E \subset L$ and all 
$ f \in \Phi_{ij},~ \mu(f(E)) = \mu (Q(E) + a) = \frac{1}{|\det Q|} \mu(E),\\
 \mbox{where} 
~{f : x \mapsto Qx + a}$. In particular, $\mu(f(W_{j})) = \frac{1}{q} ~w_{j},
~ \mbox{where}~ w_{j} := \mu (W_{j})~ \mbox{and}~  q=|\det Q| $. We obtain 
\[ w_{i} \leq \sum_{j=1}^{m} \frac{1}{q^{r}}~ \mbox{card}((\Phi^{r})_{ij})~w_{j} \]
from (\ref{10}).

Let $w = [w_{1},\cdots,w_{m}]^T$.
Since $\bigcup_{i=1}^{m} W_{i} = \overline{L}$, the Baire category theorem assures
us that for at least one $i$, 
\be
\stackrel{\circ}{W_{i}} \neq \emptyset \label{11}
\ee 
and then the primitivity gives this for all $i$.
So  $w > 0$ and
\be
 w \leq \frac{1}{q^{r}} S(\Phi^{r}) w \leq \frac{1}{q^{r}} S(\Phi)^{r} w,
~~ \mbox{for any}~ r \geq 1 \,.
\ee
Since the PF-eigenvalue of $S(\Phi)^r$ is $~q^r = |\det Q|^r$ and $S(\Phi)^r$ is primitive,  we
have from  Lemma \ref{PF1} that 
\be w = \frac{1}{q^{r}}S(\Phi^{r}) w =  \frac{1}{q^{r}}S(\Phi)^{r} w,
 ~~\mbox{for any}~ r \geq 1\, . 
\ee
The positivity of $w$ together with  $S(\Phi^{r}) \leq S(\Phi)^{r}$ shows that 
$S(\Phi^{r})= S(\Phi)^{r}$. This proves (i) and (ii).

Fix any $i \in \{1,\cdots,m\}$, let $\stackrel{\circ}{W_{i}}$ contain a basis open 
set $a + Q^{r} \overline{L}$ with some $r \in \ZZ_{\geq 0}$ by (\ref{11}).
 Since $(\tilde{U},\Phi^{r})$ is a substitution system, 
$~ a+ Q^{r}\overline{L} \subset \, \stackrel{\circ}{W_{i}} \, \subset 
W_{i} = \bigcup_{j=1}^{m} (\Phi^{r})_{ij} W_{j}$. 
In particular, $(a + Q^{r}\overline{L}) \cap g(W_{k}) \neq \emptyset\,$ for some 
$k \in \{1,\cdots,m\}$ and some $g \in (\Phi^{r})_{ik}$. 
However $~g(\overline{L}) = b + Q^{r}\overline{L}$ for some $b \in L$, 
so $(a + Q^{r}\overline{L}) \cap (b + Q^{r}\overline{L}) \neq \emptyset$. 
This means $a + Q^{r}\overline{L} = b + Q^{r}\overline{L}$. Thus 
\be 
g(W_{k}) \subset g(\overline{L}) = a + Q^{r}\overline{L} \subset \stackrel{\circ}{W_{i}}\,.\label{14}
\ee
For all $f \in (\Phi^{r})_{ij},~ j \in \{1,2,\cdots,m\}$, $f$ is clearly an open map, 
so $\bigcup_{j=1}^{m}(\Phi^{r})_{ij}(\stackrel{\circ}{W_{j}}) \subset \,\stackrel{\circ}{W_{i}}\,.$ Thus 
\begin{eqnarray}
 {\partial W_{i} = W_{i} ~\backslash \stackrel{\circ}{W_{i}}} & = &{\left( \bigcup_{j=1}^{m}(\Phi^{r})_{ij}(W_{j}) \right) ~\backslash \stackrel{\circ}{W_{i}}} \nonumber\\
& \subset &{\bigcup_{j=1}^{m} \left( (\Phi^{r})_{ij}(W_{j}) ~\backslash~ (\Phi^{r})_{ij}(\stackrel{\circ}{W_{j}}) \right)}\nonumber\\
& \subset &{\bigcup_{j=1}^{m}(\Phi^{r})_{ij}(\partial W_{j})\,.} \label{15}
\end{eqnarray} 
Note that due to (\ref{14}) at least one $g$ in $(\Phi^{r})_{ij}$ does not contribute to the relation (\ref{15}).

Let $v_{i} := \mu (\partial W_{i}),\, i = 1,\cdots,m\,$ and $v := [v_{1},\cdots,v_{m}]^T$. So $~v \leq
\frac{1}{q^{r}} S(\Phi^{r}) v$. Actually, by what we just said,
\be
0 \leq v \leq \frac{1}{q^{r}}S' v \leq \frac{1}{q^{r}}S(\Phi^{r})v = \frac{1}{q^{r}}S(\Phi)^{r}v,\label{16}
\ee 
where $S' \leq S(\Phi)^{r}, S' \neq S(\Phi)^{r}$.
Now applying the Lemma \ref{PF1} again we obtain 
equality throughout (\ref{16}). 
But by Lemma \ref{PF2} the eigenvalues of $\frac{1}{q^{r}}S'$ are strictly less in absolute value 
than the PF-eigenvalue of
$\frac{1}{q^{r}}S(\Phi)^{r}$, which is $1$.  This forces $v = 0$, and hence $\mu (\partial W_{i}) = 0, i =
1,\cdots,m$. \hfill $\square$ 

\vspace{.4cm}
In the sequel, the central concern is to relate the sets \,
$U_i$ and the sets $\Lambda_i:= W_i \cap L$. Clearly $\Lambda_i
\supset U_i$. The next 
lemma groups a circle of ideas that relate this question to the boundaries 
and interiors of the $W_i$.

\begin{lemma} \label{3lemmas}
Let $U_{i}, i = 1,\cdots,m,$ be point sets of the lattice $L$ in $\RR^{n}$.
Let $Q$ be an inflation of $L$ and identify $L$ with its image in 
its $Q$-adic completion $\overline{L}$. Define 
 $W_{i} := \overline{U_{i}}$ in $\overline{L}$ and $\Lambda_{i} := W_{i} \cap L$. 
\begin{itemize} 
\item[(i)]
If $U_1, \dots ,U_m$ are disjoint and $\mu(\overline{\Lambda_{i}
\backslash U_{i}}) = 0~~ \mbox{for all}~ i = 1,\cdots,m$, then
$\stackrel{\circ}{W_{i}} \cap \stackrel{\circ}{W_{j}} = \emptyset$  for all $ i \neq j$.

\item[(ii)]
If $L = \bigcup_{i=1}^{m} U_{i}$ and
 $\stackrel{\circ}{W_{i}} \cap \stackrel{\circ}{W_{j}} = \emptyset$  for all $ i \neq j$,
 where $ i,j \in \{1,\cdots,m\}$, then $\Lambda_{i} \backslash U_{i} \subset \bigcup_{j=1}^{m}
\partial W_{j}~~ \mbox{for all}~ i = 1,\cdots,m $.
\item[(iii)]
 If $\mu(\partial W_{j}) = 0~$ for all $j = 1,\cdots,m$ and
  $\Lambda_{i} \backslash U_{i} \subset \bigcup_{j=1}^{m} \partial W_{j}$, then
$\mu(\overline{\Lambda_{i} \backslash U_{i}}) = 0$.
\end{itemize}
\end{lemma}
{\sc{proof}}:
(i) Suppose there are $\,i,\,j\, \in \{1,\cdots,m\}$ 
with $\stackrel{\circ}{W_{i}} \cap \stackrel{\circ}{W_{j}} \neq \emptyset$. 
We can choose $a \in
{(\stackrel{\circ}{W_{i}} \cap \stackrel{\circ}{W_{j}}) \cap L}$, since $L$ is dense in $\overline{L}$ and
$\stackrel{\circ}{W_{i}} \cap \stackrel{\circ}{W_{j}}$ is open. 
Choose $k \in \ZZ_{+}$ so
that $a + q^{k}\overline{L} \subset \stackrel{\circ}{W_{i}} \cap \stackrel{\circ}{W_{j}}$. 
Note that $a +q^{k}L \subset \Lambda_{i} \cap \Lambda_{j}$. Then 
\begin{eqnarray*}
{\bigcup_{i=1}^{m} (\Lambda_{i} \backslash U_{i})}& 
\supseteq &{\left( (a + q^{k}L)~ \backslash~ U_{i} \right) 
\cup \left( (a + q^{k}L)~ \backslash~ U_{j}
\right)}\nonumber\\ 
& = &{(a + q^{k}L) ~\backslash~ (U_{i} \cap U_{j})}\nonumber\\
& = &{ a + q^{k}L,~~\mbox{since the }~U_{i},~i=1,\cdots,m,~\mbox{are disjoint}\,.}\nonumber
\end{eqnarray*}
So 
\begin{eqnarray*}
 {\sum_{i=1}^{m} \mu (\overline{\Lambda_{i} \backslash U_{i}})} & \geq & 
{\mu (\bigcup_{i=1}^{m} (\overline{\Lambda_{i} \backslash U_{i}}))}\nonumber\\
 & \geq & {\mu (a + q^{k}\overline{L})}\nonumber\\
 & > &{0\,,}\nonumber
\end{eqnarray*}
contrary to assumption. 

(ii) Assume $\stackrel{\circ}{W_{i}} \cap \stackrel{\circ}{W_{j}} = \emptyset$  
for all $i \neq j$. For any $i \in \{1,\cdots,m\}$,
\begin{eqnarray*}
 {(\Lambda_{i} \backslash U_{i})} & \subset & {( \bigcup_{j \neq i} U_{j}) \cap W_{i}, 
~~~\mbox{since}~ L = \bigcup_{i=1}^{m} U_{i} }\\
 & \subset &{\bigcup_{j \neq i} (W_{j} \cap W_{i})}
  \subset {\bigcup_{j=1}^{m} \partial W_{j}, ~~~\mbox{since}~ \stackrel{\circ}{W_{i}} \cap 
 \stackrel{\circ}{W_{j}} = \emptyset~~ \mbox{for all}~ i \neq j \,.} 
\end{eqnarray*}

(iii) Obvious. \hfill $\square$

\section{Model Sets}

Let us recall the notion of a model set (or cut and project set).
A {\em {cut and project scheme}} (CPS) consists of a collection of spaces and mappings as follows;
\be
\begin{array}{ccccc}
 \RR^{n} & \stackrel{\pi_{1}}{\longleftarrow} & \RR^{n} \times G & \stackrel{\pi_{2}}
{\longrightarrow} & G \\ 
 && \bigcup \\
 && \tilde{L}
\end{array} 
\ee
where $\RR^{n}$ is a real Euclidean space, $G$ is some locally 
compact Abelian group, and 
$ \tilde{L} \subset {\RR^{n}
\times G}$ is a lattice,  i.e.  a discrete subgroup for which the quotient group 
$(\RR^{n} \times G) / \tilde{L}$ is
compact. Furthermore, we assume that $\pi_{1}|_{ \tilde{L}}$ is injective 
and $\pi_{2}(\tilde{L})$ is dense in $G$.
   
A {\em model set} in $\RR^{n}$ is a subset of $\RR^{n}$ which, up to 
translation, is of the
 form $\Lambda(V) = \{~\pi_{1}(x)~ |~ x \in \tilde{L}, \pi_{2}(x) \in V\}$ 
for some cut and project scheme as above, where $V \subset G$ has non-empty 
interior and compact closure (relatively compact). When we need to be more
precise we explicitly mention the cut and project scheme from which
a model set arises. This is quite important in some of the theorems
below. Model sets are always Delone subsets of $\RR^n$, that is to say, 
they are relatively dense and uniformly discrete.
 
We call the model set $\Lambda(V)$ {\em regular} if the boundary  
$\partial V = \overline{V} \backslash \stackrel{\circ}{V}$ of $V$
is of (Haar) measure $0$.
We will also find it convenient to consider certain degenerate types of model
sets. A {\em weak} model set is a set in $\RR^{n}$ of the form
$\Lambda(V)$  where we assume only that $V$ is relatively compact, but
not that it has a non-empty interior. When $V$ has no interior,
$\Lambda(V)$ is not necessarily relatively dense
in $\RR^n$ but regularity still means that the boundary of $V$ is of measure $0$. 

\begin{theorem}\label{martinThm} \textmd(Schlottmann \cite{martin3}) If 
$\Lambda = \Lambda(V)$ is a regular model set, 
then $\Lambda$ is a pure point diffractive set, i.e. the Fourier transform of 
its volume averaged autocorrelation measure is a pure point measure.
\end{theorem} \hfill $\square$ 

It is this theorem that is a prime motivation for finding criteria for sets to 
be model sets.

Now let $(\tilde{U},\Phi)$ be a substitution system with inflation $Q$ on a 
lattice $L$ of $\RR^n$ and let
$\overline{L}$  be the $Q$-adic completion of $L$. This gives rise to the cut 
and project scheme. 
\be \label{cps}
\begin{array}{ccccc}
\RR^{n} & \stackrel{\pi_{1}}{\longleftarrow} & \RR^{n} \times \overline{L} & 
\stackrel{\pi_{2}}{\longrightarrow} & \overline{L}\\
&& \bigcup \\
L & \longleftarrow & \tilde{L} & \longrightarrow & L\\
t & \longleftarrow & (t,t) & \longrightarrow & t
\end{array} 
\ee 
where $\tilde{L} := \{~(t,t)~ | ~t \in L \} \subset \RR^{n} \times \overline{L}$.

We claim that $(\RR^{n} \times \overline{L}) / \tilde{L}$ is compact. $\tilde{L}$
 is clearly discrete and closed in $\RR^n \times \overline{L}$.
Since  $(\RR^{n} \times \overline{L}) / \tilde{L}$ is Hausdorff and satisfies the first axiom of
countability, it is enough to show that it is sequentially compact \cite{Kelley}.
If $\{(x_{i},z_{i}) + \tilde{L}\}$ is a countable sequence in $(\RR^{n} \times \overline{L}) / \tilde{L}$, 
then there is a subsequence $\{(x_{i},z_{i}) + \tilde{L}\}_{S}$~ with $\{x_{i} + L\}_{S}$ 
convergent sequence,  since $\RR^{n} / L$ is compact. We can rewrite 
$\{(x_{i},z_{i}) + \tilde{L}\}_{S}$ as $\{(x_{i}',z_{i}') + \tilde{L}\}_{S}$, 
where $\{x_{i}'\}_{i \in S}$ converges to $x$ in $\RR^{n}$. 
Since $\overline{L}$ is compact, there is a convergent subsequence $\{z_{i}'\}_{S'}$ to 
some $z$ in $\overline{L}$. 
Thus $\{(x_{i}',z_{i}')\}_{S'}$ converges to $(x,z)$ in $\RR^{n} \times \overline{L}$. 
Therefore $(\RR^{n} \times \overline{L}) / \tilde{L}$ is sequentially compact. 

Note also that $\pi_{1}|_{\tilde{L}}$ is 
injective and $\pi_{2}(\tilde{L})$ is dense in $\overline{L}$. 

\begin{lemma} \label{regModelSet}
Let $U_{i}, i = 1,\cdots,m,$ be disjoint point sets of the lattice $L$ in 
$\RR^{n}$. Identify $L$ and its image in $\overline{L}$. 
Let $W_{i} := \overline{U_{i}}$ in $\overline{L}$ and $\Lambda_{i} := W_{i} \cap L$.
 Suppose that $\mu(\partial W_{i}) = 0$ for all $i = 1,\cdots,m$. 

(i) If $\Lambda_{i} \backslash U_{i} \subset \bigcup_{j=1}^{m} \partial W_{j}$ 
then, relative to the CPS({\rm \ref{cps}}), $U_i$ is a regular weak model set when 
$\stackrel{\circ}{W_{i}}$ is empty, and $U_{i}$ is a regular model set when
$\stackrel{\circ}{W_{i}}$ is non-empty.

(ii) If  $L = \bigcup_{j=1}^{m} U_{j}$ and each $U_i$ is a regular model set,
then  $\Lambda_{i} \backslash U_{i} \subset \bigcup_{j=1}^{m} \partial W_{j}$
for all $i=1,\dots,m$.
\end{lemma}
{\sc {proof}}: (i) 
Assume that $\Lambda_{i} \backslash U_{i} \subset \bigcup_{j=1}^{m} \partial W_{j}$ for all 
$i = 1,\cdots,m$.
Since $\mu(\partial W_{i}) = 0$ for all $i = 1,\cdots,m$, 
\be \label{measW}
\mu(W_{i}) = \mu(\stackrel{\circ}{W_{i}}) = 
\mu(\stackrel{\circ}{W_{i}} \backslash \bigcup_{j=1}^{m} \partial W_{j}) \label{51}
\ee

Since $\Lambda_{i} = W_{i} \cap L$, $U_{i} = V_{i} \cap L$ 
where $V_{i} := W_{i} \backslash (\Lambda_{i}
\backslash U_{i})$. 
Now $V_{i} \supset\, \stackrel{\circ}{W_{i}} \backslash \bigcup_{j=1}^{m} \partial W_{j}$.
From $\stackrel{\circ}{W_{i}} \backslash \bigcup_{j=1}^{m} \partial W_{j} \subset \,
 \stackrel{\circ}{V_{i}} \subset V_{i} \subset \overline{V_{i}} = W_{i}$ 
and (\ref{measW}), $\mu(\overline{V_{i}} \backslash
\stackrel{\circ}{V_{i}}) = 0$. So $U_i$ is regular. 
If $\stackrel{\circ}{W_{i}} = \emptyset$, then $\stackrel{\circ}{V_{i}} = 
\emptyset$ also. Thus $U_i$ is a regular weak model set.
On the other hand, for any $i$ with $\stackrel{\circ}{W_{i}} \neq \emptyset$, 
$\stackrel{\circ}{V_{i}} \neq \emptyset$ and 
$\overline{V_{i}}$ is compact. It follows that $U_{i} = \Lambda(V_{i})$ is a regular model set
for the  CPS (\ref{cps}).

(ii) Suppose that  $\stackrel{\circ}{V_{i}} \neq \emptyset$,
$\mu(\overline{V_{i}} \backslash \stackrel{\circ}{V_{i}}) = 0$, 
where $U_{i} = V_{i} \cap L$, and $L = \bigcup_{j=1}^{m} U_{j}$. 
Then from $\overline{\Lambda_{i} \backslash U_{i}} =
\overline{\Lambda(W_{i})
\backslash \Lambda(V_{i})} \subset \overline{W_{i} \backslash V_{i}} \subset W_{i} \backslash
\stackrel{\circ}{V_{i}} = \overline{V_{i}} \backslash \stackrel{\circ}{V_{i}}$, we have
$\mu(\overline{\Lambda_{i} \backslash U_{i}}) = 0$ for all $i = 1,\cdots,m$. By 
Lemma \ref{3lemmas} (i) and (ii) ,
$\stackrel{\circ}{W_{i}} \cap \stackrel{\circ}{W_{j}} = 0$ for all $i
\neq j$ and  $\Lambda_{i} \backslash U_{i} \subset \bigcup_{j=1}^{m} \partial W_{j}$. \hfill $\square$

\begin{theorem} \label{mainTheorem}
  Let $(\tilde{U},\Phi)$ be a primitive substitution system with inflation $Q$ on the lattice $L$ 
in $\RR^{n}$. Suppose that PF-eigenvalue of the substitution matrix $S(\Phi)$ 
is equal to $|\det Q|$ and $L = \bigcup_{i=1}^{m} U_{i}$. Then the following are equivalent. 

(i) There is a primitive substitution matrix $\Psi$ admitting a coincidence,  
where $(\tilde{U},\Psi)$ is equivalent to $(\tilde{U},\Phi^{M}) $ for some $M \geq 1$.

(ii) The sets $U_{i},\; i=1,\cdots,m,$ of $\tilde{U}$ are model sets for  the CPS (\ref{cps}).

(iii) For at least one $i$, $U_{i}$ contains a coset \, $a + Q^{M}L$.

(iv) $(\tilde{U},\Phi)$ admits a modular coincidence.
\end{theorem}
{\sc{proof}} : 

(i) $\Rightarrow$ (ii): 
Suppose that $(\tilde{U},\Psi)$ admits a coincidence and is equivalent to $(\tilde{U},\Phi^{M})$.
Fix $i \in \{1,\cdots,m\}$ with $\cap_{j=1}^{m} \Psi_{ij} \neq \emptyset$ and let $g$ be in this 
intersection. Recalling equation (\ref{10}), and in view of the choice of $g$, we have
\[ \mu(W_{i}) \leq \left( \sum_{j=1}^{m} \sum_{f \in \Psi_{ij}} \mu (f(W_{j}))\right) - \mu (g(W_{k}) 
\cap g(W_{l}))\,,\]
for any $k,l \in \{1,\cdots,m\}$ with $k \neq l$.
On the other hand, from Theorem \ref{mainPFtheorem} (ii)
\be \mu(W_{i}) = \frac{1}{q^{M}} \sum_{j=1}^{m} (S(\Psi))_{ij} \mu (W_{j}) = 
\sum_{j=1}^{m} \sum_{f \in \Psi_{ij}} \mu (f(W_{j}))\,.
\ee
Thus, in fact, $\mu (g(W_{k}) \cap g(W_{l})) = 0$  whenever $k \neq l$. 
It follows at once that $\stackrel{\circ}{W_{k}} \cap \stackrel{\circ}{W_{l}} = \emptyset ~$ 
for all $k \neq l$,
 since the measure of any open set is larger than $0$. 

Recall that $\stackrel{\circ}{W_{i}} \neq
\emptyset$  and $\mu(\partial W_{i}) = 0$ for all $i = 1,\cdots,m$ . Then by 
Lemma \ref{3lemmas}(ii) and Lemma \ref{regModelSet},
$U_{i}, \,  i = 1,\cdots,m,$ are model sets in CPS(\ref{cps}).
 
(ii) $\Rightarrow$ (iii): 
Assume that  $U_{i}, i = 1,\cdots,m,$ are model sets in CPS(\ref{cps}), i.e. $U_{i} = \Lambda(V_{i}) = V_{i}
\cap L$ for some $V_{i}$ with $\stackrel{\circ}{V_{i}} \neq \emptyset$. Thus there is a coset $a +
Q^{M}\overline{L} \subset \stackrel{\circ}{V_{i}}$ and, since we can always choose the coset
representative from the dense lattice $L$, we can arrange that $a + Q^{M}L \subset U_{i}$. 

(iii) $\Rightarrow$ (iv):
Assume that for at least one $i$, $U_{i}$ contains a coset $a + Q^{M}L$. Fix $i$.
 Iterate $\Phi$ M-times. Then each function $f$ in the substitution system $\Phi^{M}$ has the 
form $f : x \mapsto Q^{M}x + b$. 
For each $j$, let $G_{j} := \{~f \in (\Phi^{M})_{ij}~ |~ t(f) \equiv a ~~
\mbox{mod}~ Q^{M}L\}$. (Recall that $t(f)$ is the translational part of $f$).
From $U_{i} = \bigcup_{j=1}^{m} \bigcup_{f \in (\Phi^{M})_{ij}} f(U_{j})$, we
obtain $a + Q^{M}L \subset \bigcup_{j=1}^{m} \bigcup_{f \in G_{j}}f(U_{j})$. In fact 
\be
a + Q^{M}L = \bigcup_{j=1}^{m} \bigcup_{f \in G_{j}}f(U_{j})\,,
\ee
since the right hand side is clearly inside $a + Q^{M}L$. 
From the fact $a + Q^{M}L \subset U_{i}$, we get $
\Phi^{M}[a] = \bigcup_{j=1}^{m} G_{j} \subset \bigcup_{j=1}^{m} (\Phi^{M})_{ij}$. 
Therefore $\Phi^{M}$ has
a row containing an entire congruence class $\Phi^{M}[a]$. 

(iv) $\Rightarrow$ (i):
Assume $\Phi^{M}$ has a row, say $i$-th row, containing an entire congruence class $\Phi^{M}[a]$. Let
$G_{j} := \Phi^{M}[a] \cap (\Phi^{M})_{ij}$. Then 
$ \bigcup_{j=1}^{m} \bigcup_{f \in G_{j}}f(U_{j}) \subset a + Q^{M}L$.
Recall that $\bigcup_{j=1}^{m} U_{j} = L$ and $\tilde{U} = \Phi^{M}(\tilde{U})$. It follows that 
the elements of $a + Q^{M}L$ can be obtained from the substitution system $\Phi^{M}$ {\em only} 
from the mappings of $\Phi^{M}[a]$, and indeed they must {\em all} appear as images of the mappings of
$\Phi^{M}[a]$. Thus 
\be
a + Q^{M}L = \bigcup_{j=1}^{m} \bigcup_{f \in G_{j}}f(U_{j}) \subset U_{i} \,. \label{23}
\ee 
On the other hand, 
\be
a + Q^{M}L = \bigcup_{j=1}^{m} Q^{M}(U_{j}) + a \,, \label{24}
\ee 
which is a disjoint union.  

We now alter our substitution system $\Phi^{M}$ as follows: 
Define $g : L \rightarrow L$ by $g(x) = Q^{M}x + a$. We may, by restriction of domain, consider $g$ 
as a function on $U_{j}, ~j = 1,\cdots,m$. We define $\Psi$ by
\[ 
 \left\{   \begin{array}{l}
           \Psi_{ij} = ((\Phi^{M})_{ij} ~\backslash~ G_{j}) \bigcup \{g\}\\
           \Psi_{kj} = (\Phi^{M})_{kj} ~~~ \mbox{if} ~k \neq i \,,
           \end{array} 
 \right.  \]
for all $j$.
From (\ref{23}) and (\ref{24}), the $\Psi_{ij}, j = 1,\cdots,m$, consist of maps from $U_{j}$ to $U_{i}$
 and have the same total effect on $U_{i}$ as the $(\Phi^{M})_{ij},j = 1,\cdots,m\,$. 
 Thus $(\tilde{U},\Psi)$ is a substitution system admitting a coincidence.

 Since $S(\Phi^{M})$ is primitive, the incidence matrix $I(\Phi^{M})$ is primitive.  
Then $I(\Psi)$ is also primitive, since $I(\Phi^{M}) \leq I(\Psi)$. So $\Psi$ is primitive.
In addition,  $\Psi$ has the inflation $Q^{M}$ for $L$ which is an inflation in $\Phi^{M}$. 

 We claim that $S(\Psi), S(\Phi^M)$ have the same PF-eigenvalue and right 
PF-eigenvector. Then $(\tilde{U},\Psi)$ is equivalent to $(\tilde{U},\Phi^{M})$.

 We verify first that 
$\stackrel{\circ}{W_{k}} \cap \stackrel{\circ}{W_{j}} = \emptyset ~\mbox{for all}~ k \neq j$.
We can assume that $m > 1$, since there is nothing to prove when $m = 1$. 
Let $g_1 \in G_{l} = (\Phi^{M})_{il}[a] \neq \emptyset\,$ for some $l$. Take any $k \in \{1,\cdots,m\}$. 
There is $M_{0} \in \ZZ_{+}$ for which $(\Phi^{M_{0}})_{lk} \neq \emptyset$. Choose 
$f \in (\Phi^{M_0})_{lk}$.
Let $g_1 : x \mapsto Q^{M}x + a_{1}$, where $a_{1} \equiv a ~\mbox{mod}~ Q^{M}L$, and 
$f : x \mapsto Q^{M_{0}}x + b\,$ with $b \in L$.
Then $g_1 \circ f : x \mapsto Q^{M+M_{0}}x + Q^{M}b + a_{1}$. 
So $g_1 \circ f \in (\Phi^{M+M_{0}})_{ik}[a_{1} + Q^{M}b]$.
Furthermore $\,(a_{1} + Q^{M}b) + Q^{M+M_{0}}(L) \subset a_{1}+Q^{M}L \subset U_{i}$.

Let $N:= M+M_{0},\, c := a_{1}+Q^{M}b,\, \mbox{and}\, p:=g_1 \circ f$.
Note that 
\be \label{unionEqn}
c + Q^{N}L = \bigcup^{m}_{j=1} \bigcup_{h \in H_{j}} h(U_{j})\, ,
\ee 
where $H_{j} = (\Phi^{N})_{ij}[c].$

There are at least two functions in $\bigcup^{m}_{j=1}H_{j}$, since 
for all $j$~ $U_{j}, \neq L$. 
We can write  $c + Q^{N}L$ in the form
\be \label{basicDecomp}
c + Q^{N}L  = \bigcup \{~Q^{N}U_{j} + Q^{N}\alpha_{h} + c~ |~ j \in \{1,\cdots,m\}, h \in H_{j},
\alpha_{h} \in L~\}\, ,
\ee   
where we have used the explicit form of each of the mappings $h \in H_j$. This union is 
disjoint, and as a consequence the elements $\alpha_h \in L$ for $h$ in any single $H_j$
are all distinct. In particular we have $\alpha_p$ coming from $ H_k$. From (\ref{basicDecomp})
we have 
\be
L  = \bigcup_{j=1}^m \bigcup_{h \in H_j} (U_{j} + \alpha_{h} )
\ee
and separating off $U_k$,
\be \label{separating}
L  =  U_k \cup \bigcup_{j=1}^m \bigcup_{h \in H'_j} (U_{j} + \alpha_{h} -\alpha_p ) \, ,
\ee
where $H'_j := H_j$ if $j \ne k$ and $H'_k := H_k\backslash \{p\}$. Again these decompositions
are disjoint. But we also know that $U_{k}$ and $\bigcup^{m}_{\stackrel{j=1}{j \neq k}} U_{j}$ are
disjoint, and it follows that
\[ \bigcup^{m}_{\stackrel{j=1}{j \neq k}} U_{j} \subset \bigcup^{m}_{j=1} 
\bigcup_{h \in H'_{j} }(U_{j} + \alpha_{h} - \alpha_{p}).
\]
Taking closures,
\be \label{closures}
\bigcup^{m}_{\stackrel{j=1}{j \neq k}} W_{j} \subset \bigcup^{m}_{j=1} \bigcup_{h 
\in H'_{j}}(W_{j} +
\alpha_{h} - \alpha_{p}).
\ee
On the other hand, if we apply Theorem \ref{mainPFtheorem}(ii) to $\Phi^N$ and look
at (\ref{unionEqn}) we see that 
\[
\mu(c + Q^{N}\overline{L}) = \sum^{m}_{j=1}\sum_{h \in H_{j}} \mu(h(W_{j})) = \sum^{m}_{j=1}\sum_{h \in H_{j}} \mu(Q^{N}(W_{j} + \alpha_{h}) + c),
\]  and hence
\[
\mu(\overline{L}) = \sum^{m}_{j=1}\sum_{h \in H_{j}} \mu (W_{j} + \alpha_{h}) = 
\sum^{m}_{j=1}\sum_{h \in H_{j}} \mu(W_{j} + \alpha_{h} - \alpha_{p}).
\]
 Thus 
\[\mu(\overline{L}) = \mu(W_{k}) +  \left( \sum^{m}_{j=1} \sum_{h \in H'_{j} } \mu(W_{j} +
\alpha_{h} - \alpha_{p}) \right) 
\]
which, after taking closures in  (\ref{separating}),  gives
\be
\mu \left( W_{k} \cap \left( \bigcup^{m}_{j=1} \bigcup_{h \in H'_{j} }(W_{j} + \alpha_{h} - \alpha_{p})
\right) \right) = 0\,. \label{25}
\ee 
Finally from (\ref{closures}) and (\ref{25}) we obtain
\[ \mu(W_{k} \cap (\bigcup^{m}_{\stackrel{j=1}{j \neq k}}W_{j})) = 0,
\]
from which $\stackrel{\circ}{W_{k}} \cap \stackrel{\circ}{W_{j}} = \emptyset ~$ for all $k \neq j$.
 This establishes the claim.
  
Now
\begin{eqnarray} 
{\mu \left( \bigcup_{j=1}^{m} g(W_{j}) \right) }& = &
{\frac{1}{|\det Q^{M}|} \mu \left( \bigcup_{j=1}^{m} W_{j} \right) }\nonumber\\
& = &{\frac{1}{|\det Q^{M}|} \sum^{m}_{j=1}\mu(W_{j})\,,}\nonumber\\
&&{\mbox{from}~ \mu(\partial W_{j})=0,
~\stackrel{\circ}{W_{i}} \cap \stackrel{\circ}{W_{j}} = \emptyset~ \mbox{for all}~ i \neq j }\nonumber\\
& = &{\sum^{m}_{j=1} \mu(g(W_{j}))}\,. \label{26}
\end{eqnarray}
Again using Theorem \ref{mainPFtheorem} (ii), this time for $\Phi^M$, we obtain
\[
~ w = \frac{1}{|\mbox{det} Q^{M}|} S(\Phi^{M}) w \, ,
\] 
where $w = [\mu(W_1), \dots, \mu(W_m)]^T$. 
The part of this relation  in $W_i$ which pertains to the coset $a+Q^M\overline{L} $ is
\be \mu (a + Q^{M}\overline{L}) = \sum_{j=1}^{m} \sum_{f \in G_{j}} \mu (f(W_{j}))\,.
\label{27.1}
\ee 
But from (\ref{24})
\be
\mu (a + Q^{M}\overline{L}) = \mu \left( \bigcup_{j=1}^{m} g(W_{j}) \right)\,. 
\label{27.2}
\ee
Together, (\ref{26}), (\ref{27.1}), and (\ref{27.2}) show 
\[
w = \frac{1}{|\mbox{det}Q^{M}|} S(\Psi) w \,. 
\]
Since $w > 0$ and $S(\Psi)$ is primitive, $S(\Psi)$ has PF-eigenvalue 
$|\mbox{det} Q^{M}|$ and PF-eigenvector $w$ as required. \hfill $\square$

\vspace{.4cm}
{\bf Remark}: Let $A = \{a_1, \dots, a_m \}$ be an alphabet of $m$
symbols and let $\sigma$ be a primitive equal-length alphabetic
substitution system on $A$, that is,
\begin{itemize}
\item[(i)] $\sigma : A \longrightarrow A^q$ for some $q \in \ZZ_+$;
\item[(ii)] the $m \times m$ matrix $S = (S_{ij})$, whose $i,j$ entry 
is the number of appearances of $a_i$ in $\sigma(a_j)$, is
primitive.
\end{itemize}

According to Gottschalk \cite{gottschalk}, for some iteration
$\sigma^k$ of $\sigma$, there is a word $w \in A^\ZZ$ which
is fixed by $\sigma$ in the sense that 
\be
\begin{array}{ccc}
\sigma^k(w_0w_1 \dots ) &=& w_0w_1 \dots \\
\sigma^k(\dots  w_{-2} w_{-1}) &=& \dots w_{-2} w_{-1} \, .
\end{array}
\ee
Replacing $\sigma^k$ by $\sigma$ and $q^k$ by $q$ if necessary
we can suppose that $k=1$, and assume then that $\sigma(w) = w$.

We can view $w$ as a tiling of $\RR$ by tiles of types
$a_1, \dots a_m$, all of the same length $1$.
If we coordinatize each tile by its lefthand end point
so that $w_l$ gets coordinate $l$, then we obtain a partition
$U_1 \cup \dots \cup U_m$ of $\ZZ$ and an $m \times m$ matrix substitution system
$\Phi$ of $q$-affine mappings derived directly from $\sigma$:
namely, $\sigma a_j = a_{i_1} \dots a_{i_q}$ gives rise to the
mappings $(x\mapsto qx+l - 1) \in \Phi_{i_l j}$,  $l = 1, \dots, q$.

We take as our cut and project scheme

\be \label{elscps}
\begin{array}{ccccc}
\RR & \longleftarrow & \RR \times \ZZ_q & \longrightarrow & \ZZ_q\\
&& \bigcup \\
\ZZ & \longleftarrow & \tilde{\ZZ} & \longrightarrow & \ZZ\\
z & \longleftarrow & (z,z) & \longrightarrow & z
\end{array} 
\ee 
(see \ref{cps}), where $\ZZ_q$ is the $q$-adic completion of $\ZZ$.

According to Theorem \ref{mainTheorem}, the $U_i$ are model
sets for (\ref{elscps}) if and only if for some iteration
$\sigma^M$ of $\sigma$, there is a $k \in \ZZ$ for which all the 
mappings $f_l: x \mapsto q^M x + l$ with $l \equiv k \mod{q^M}$
lie in one row of $\Phi^M$.

Since $\sigma^M a_j$ has $q^M$ letters in it, there are
$q^M$ mappings in the $j$th column of $\Phi^M$.
Furthermore, since the letters $\sigma^M a_j$ are represented by 
contiguous tiles, their coordinates fall in a range of consecutive
integers, and so the mappings of the $j$th column of $\Phi^M$ are
the maps $f_l$, where  $0 \le l < q^M$, in some order. In particular,
all of the mappings in $\Phi^M$ are of this restricted form. It follows
that modular coincidence is equivalent to the existence of a row
of $\Phi^M$, say the $i$th row, and a $k$, $0 \le k < q^M$, so
that $f_k$ belongs to {\em each} of $\Phi^M_{i1}, \dots, \Phi^M_{im}$.

This condition precisely says that there is a $k$ so that the 
$k$th position of $\sigma^M(a_j)$ contains the same letter $a_i$
for all $j$. This is the well-known coincidence condition of
Dekking \cite{dekking}, and he has proved that for non-periodic
primitive equal-length substitutions, this condition is equivalent
to pure point diffractivity. It is straighforward to show
that $S(\Phi)$ has its PF-eigenvalue equal to $|\det Q|$. 
Thus we have 

\begin{coro}
Let $\sigma$ be a primitive equal-length ($= q$) alphabetic
substitution with a fixed bi-infinite word $w$, and assume
that $w$ is not periodic. Let $\Phi$ be the corresponding
matrix substitution system and let $\ZZ = U_1 \cup \dots \cup U_m$
be the corresponding partition of $\ZZ$. Then the following 
are equivalent:
\begin{itemize}
\item[(i)] there is an $M$ so that $\sigma^M$ has a coincidence
in the sense of Dekking;
\item[(ii)] $\Phi$ has a modular coincidence;
\item[(iii)] the $U_i$'s are model sets for (\ref{elscps});
\item[(iv)] the $U_i$'s are pure point diffractive.
\end{itemize}
\end{coro}

We note that this interesting {\em equivalence} of model sets
and pure point diffractivity is more than we can yet prove
in the higher dimensional substitution systems.

\section{Sphinx tiling}
{\footnotesize
\flushright{
Long we sought the wayward lynx\\
And bowed before the subtle sphinx\\
But solved we not the cryptic sphinx\\
Before we found the wayward links.\\
-Anon\\
}}
            
\begin{figure}[ht]
\centerline{\epsfysize=30mm \epsfbox{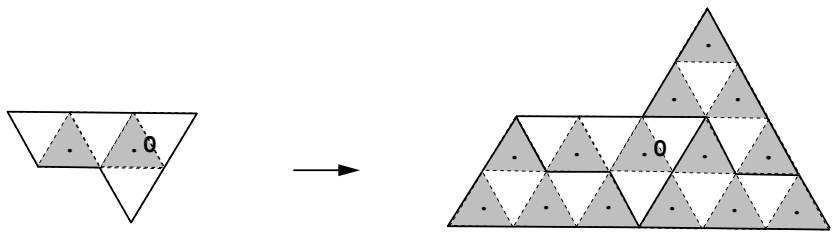}}
\caption{Sphinx Inflation [Type 1]}  \label{sphinx:1}
\end{figure}

\begin{figure}[ht]
\centerline{\epsfysize=30mm \epsfbox{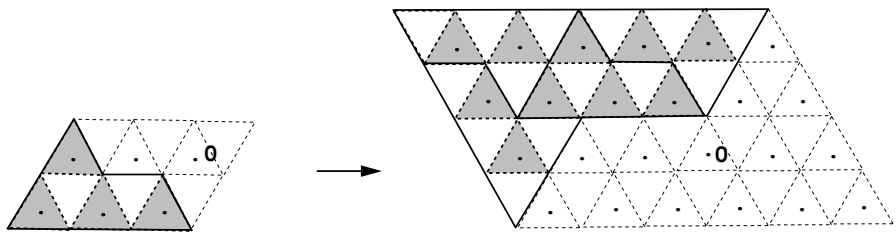}}
\caption{Sphinx Inflation [Type 2]}  \label{sphinx:2}
\end{figure}

In this section we take up the sphinx tiling. This is a substitution tiling whose 
subdivision rule is shown in Figure ~\ref{sphinx:1} and Figure ~\ref{sphinx:2}. 

It has $12$ sphinx-like tiles (up to translation).
If we choose a single point in the same way in each sphinx then we arrive at 
$12$ sets of points. We wish to show that each of these sets is a regular
model set. Actually we make a slight alteration to this, choosing several
points from each tile, but this is equivalent to our original problem.  

\begin{figure}[ht]
\centerline{\epsfysize=100mm \epsfbox{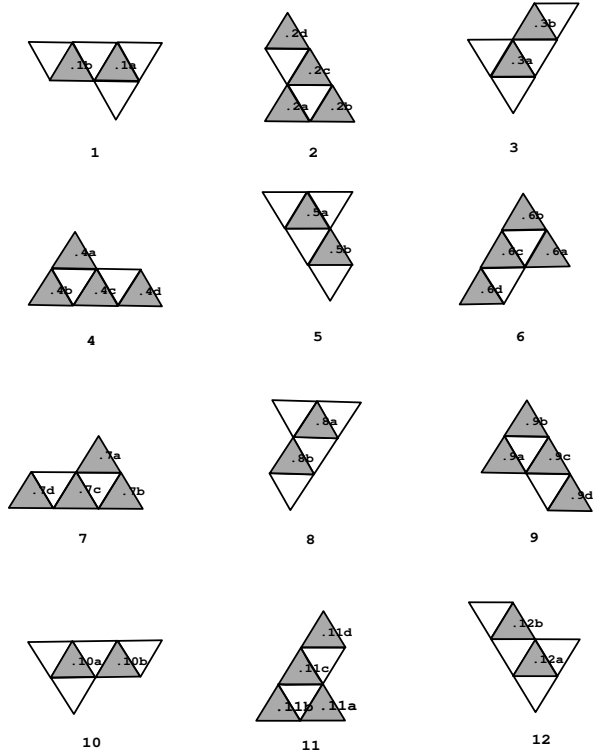}}
\caption{12 Sphinx Tiles}  \label{sphinx:3}
\end{figure}

Each sphinx can be viewed as consisting of $6$ 
equilateral triangles of two orientations. In this way, any sphinx tiling
determines a tessellation of the plane by equilateral triangles.
We consider the centre points of the triangles of one orientation.
These clearly form a lattice $L$, once we have chosen one of them
as the origin. Note that some sphinxes have two points and others have 
four points in $L$. We give names to each tile and the points in it as 
shown in Figure~\ref{sphinx:3}. Then the $12$ types of sphinx partition $L$
into 36 subsets forming a matrix substitution system. We
show that these are model sets for a $2$-adic-like cut and
project scheme of the form of (\ref{cps}).

With the origin as shown, the coordinates are chosen so that 
in the standard rectangular system $(1,0)$ is the lattice point
directly to the right of $(0,0)$. It is more convenient to
replace this by an oblique coordinate system:
$L = \{~ae + bw~ |~ a,b \in \ZZ \}$, 
where $e = (1,0), w = (\frac{1}{2},\frac{\sqrt{3}}{2})$ in the standard
rectangular system and 
relative to this basis we can identify $L$ and $\ZZ^2$ 
and denote $ae + bw$ by $(a,b)$.
The basic inflation shown in Figure~\ref{sphinx:1} gives rise
to the map 
\[T : x \mapsto 2Rx + (1,0)\, , \] 
where $R$ is a reflection in $\RR^{2}$ through $x$-axis, i.e. in the new coordinates, 
$~R(1,0) = (1,0), R(0,1) = (1,-1)$. 

The  various types of points are designated by letter pairs
$i\alpha$, where $i \in \{1,\cdots,12\}$
and $\alpha \in \{a,\cdots,d\}~$ (of which only 36 actually occur).
Let $U_{i\alpha}$ be the set of points of type $i\alpha$. 
On the basis of this we can make mappings of each point set to other point set.

Define
\[
\begin{array}{ll}
h_{1} : x \mapsto Tx + (0,0), & h_{2} : x \mapsto Tx + (1,0)\\
h_{3} : x \mapsto Tx + (0,1), & h_{4} : x \mapsto Tx + (-1,1)\\
h_{5} : x \mapsto Tx + (-1,0), & h_{6} : x \mapsto Tx + (0,-1)\\
h_{7} : x \mapsto Tx + (1,-1), & h_{8} : x \mapsto Tx + (2,-1)\\
h_{9} : x \mapsto Tx + (-1,2), & h_{10} : x \mapsto Tx + (-1,-1).
\end{array}
\]
Let $f_{i\alpha\,j\beta}$ be the function which maps $j\beta$-point set into $i\alpha$-point set.

[Type 1]
\[
\begin{array}{ll} 
f_{9a\,1a} = h_{4} : x \mapsto Tx + (-1,1), & f_{1a\,1b} = h_{2} : x \mapsto Tx + (1,0)\\ 
f_{9b\,1a} = h_{9} : x \mapsto Tx + (-1,2), & f_{1b\,1b} = h_{1} : x \mapsto Tx + (0,0)\\
f_{9c\,1a} = h_{3} : x \mapsto Tx + (0,1), &  f_{4a\,1b} = h_{5} : x \mapsto Tx + (-1,0)\\
f_{9d\,1a} = h_{2} : x \mapsto Tx + (1,0), &  f_{4b\,1b} = h_{10} : x \mapsto Tx + (-1,-1)\\
f_{4a\,1a} = h_{1} : x \mapsto Tx + (0,0), & f_{4c\,1b} = h_{6} : x \mapsto Tx + (0,-1)\\
f_{4b\,1a} = h_{6} : x \mapsto Tx + (0,-1), & f_{4d\,1b} = h_{7} : x \mapsto Tx + (1,-1)\\
f_{4c\,1a} = h_{7} : x \mapsto Tx + (1,-1), &\\
f_{4d\,1a} = h_{8} : x \mapsto Tx + (2,-1). &\\
\end{array}
\] 

[Type 2]
\[
\begin{array}{ll}
f_{12a\,4a} = h_{1} : x \mapsto Tx + (0,0), & f_{1a\,4b} = h_{2} : x \mapsto Tx + (1,0)\\
f_{12b\,4a} = h_{4} : x \mapsto Tx + (-1,1), & f_{1b\,4b} = h_{1} : x \mapsto Tx + (0,0)\\ \\
f_{4a\,4c} = h_{1} : x \mapsto Tx + (0,0), & f_{1a\,4d} = h_{1} : x \mapsto Tx + (0,0)\\
f_{4b\,4c} = h_{6} : x \mapsto Tx + (0,-1), & f_{1b\,4d} = h_{5} : x \mapsto Tx + (-1,0)\\
f_{4c\,4c} = h_{7} : x \mapsto Tx + (1,-1), &\\
f_{4d\,4c} = h_{8} : x \mapsto Tx + (2,-1). &\\
\end{array}
\]

All points in a sphinx having 2-points in it are mapped as in [Type 1] 
changing the translation part 
according to the orientation of the sphinx relative to sphinx 1. 
Likewise, all points in 
a sphinx having 4-points in it are mapped as in [Type 2] relative to sphinx 4.

\begin{figure}[ht] \label{sphinxMatrix}
\centerline{\epsfysize=80mm \epsfbox{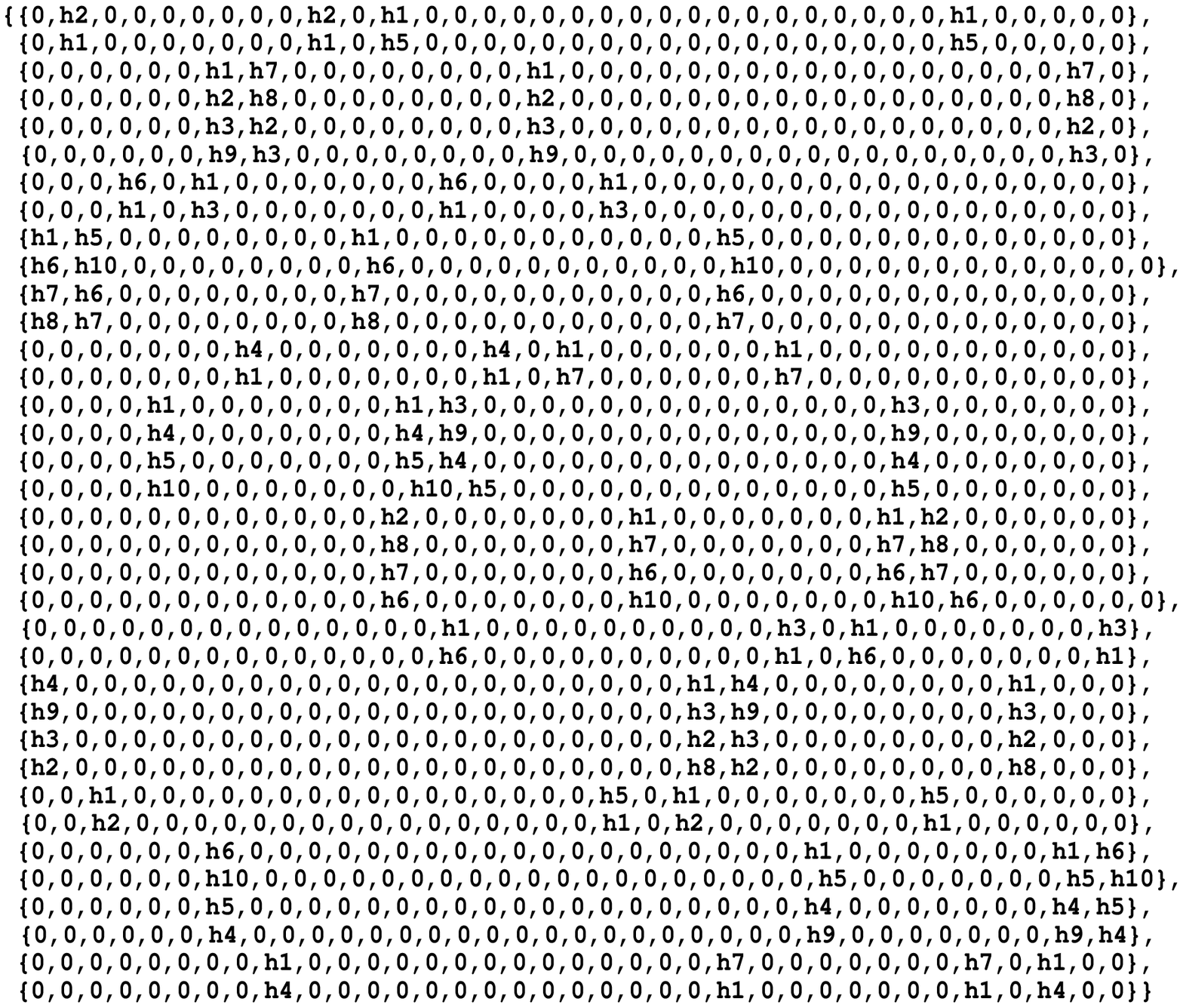}}
\caption{Sphinx matrix function system($\Phi$)}  \label{sphinxmat}
\end{figure}

Now we can list the $36 \times 36$ matrix($\Phi$) of affine mappings that make up 
our substitution system (Figure \ref{sphinxmat}).

We can check that $S(\Phi)$ has PF-eigenvalue $4$ and is a primitive matrix 
and the union of point sets is L. We used Mathematica to check
that property (iv) in Theorem~\ref{mainTheorem} is satisfied 
in $\Phi^8$ (it may actually be satisfied at some lower power).
Certainly in $\Phi^8$ there are a large number of modular coincidences.
Theorem~\ref{mainPFtheorem} and \ref{mainTheorem} say that all 36 point sets are 
regular model sets in CPS (\ref{cps}).

\section{The total index and model sets}

In this section we derive another criterion for determining when a partition 
of a lattice is a partition into $Q$-adic model sets, the difference this time being
that there is no substitution system involved. 

We assume that we are given a lattice $L$ in $\RR^n$ and an inflation $Q$ on $L$
as in (\ref{inflation}). The notation remains the same as before.  The main ingredient
is a non-negative sub-additive function called the total index which is defined
on the subsets of $L$ and its $Q$-adic completion ${\overline L}$.

For any subset $V$ of $L$ the {\em coset part} of $V$ is defined as 
\be
\mathcal{C}(V) := 
\bigcup \{\,C\,|\,C\, \mbox{is a coset in}\, V \,\} \,.
\ee
The key point to remember in what follows is that two cosets in 
$L\,(\overline{L})$ are
either disjoint or one of them is contained in the other. If $C = a +Q^kL$
is a coset then we write $[L:C]$ for the index of the subgroup $Q^kL$ in $L$.

\begin{lemma} \label{cosetDecomp}
The coset part of $V$ can be written as a disjoint union of cosets in $V$.
\end{lemma}
{\sc {proof}}:
 If $V$ contains no cosets, then the result is clear. Suppose $V$ contains cosets. 
Let $C_{1} = a_{1} + Q^{k_{1}}L$ be a coset in $V$ with $k_{1}$ minimal.
 Consider $V \backslash C_{1}$.  No coset can be partly in $C_1$ and partly 
in $V \backslash C_{1}$. Thus, if $V \backslash C_{1}$
contains no cosets,  then $\mathcal{C}(V) = C_{1}$. 
Otherwise let $C_{2}$ be a coset $a_{2} + Q^{k_{2}}L$ with $k_{2}$ minimal in 
$V \backslash C_{1}$. Then 
$\mathcal{C}(V) \supset C_{1} \stackrel{\cdot}{\bigcup} C_{2}$.
 We continue this process. Since there are only finitely many cosets 
for $Q^{k}L$ in $L$, 
either we obtain 
$\mathcal{C}(V) = C_{1} \stackrel{\cdot}{\bigcup} \cdots \stackrel{\cdot}{\bigcup}
C_{r}$ for some $r$ or 
$\mathcal{C}(V) \supset C_{1} \stackrel{\cdot}{\bigcup} C_{2}
\stackrel{\cdot}{\bigcup} \cdots,$ where $k_{1} \leq k_{2} \leq \cdots$ 
is infinite and unbounded. In the latter case, 
$\mathcal{C}(V) = \bigcup_{i=1}^{\infty} C_{i}$ is our required decompostion. 
If not, there is a  coset $C =
a + Q^{k}L$ in $V$ such that $C \not\subset \bigcup_{i=1}^{\infty}C_{i}$.
Then there is $C_{i}$ with $k_{i-1} \leq k < k_{i}.$ 
This contradicts the choice of $C_{i}.$  \hfill $\square$
\\

For $V \subset L$, we call a decomposition $\mathcal{C}(V) = \bigcup_i C_{i}$
of $\mathcal{C}(V)$ into mutually disjoint cosets using the 
algorithm of Lemma \ref{cosetDecomp},
an {\em efficient} decomposition of $V$ into cosets.
In this case we call $c(V) := \sum_i [L : C_{i}]^{-1}$ the {\em total index}
of $V$. Since any coset is an efficient decomposition of
itself, we have $c(V) = \sum_i c(C_i)$. We will see shortly that the 
total index is finite.

It is useful to note that an efficient decomposition of 
$\mathcal{C}(V) = \bigcup_i C_{i}$ of $\mathcal{C}(V)$ into cosets has
the following special property: if $D$ is any coset of $V$ then
necessarily $D \subset C_i$ for some $i$.

\begin{lemma}
Any two efficient decompositions of $\mathcal{C}(V)$ are the same up to
rearrangement of the order of the cosets. In particular the total
index is well-defined. 
\end{lemma}
{\sc proof}:
Let $\mathcal{C}(V) = \bigcup C'$ 
be a second decomposition of $\mathcal{C}(V)$ determined by the same algorithm as 
in Lemma \ref{cosetDecomp}. Then with $k_{1}$ as 
in the Lemma, let $D_{1},\cdots,D_{r}$ be all the cosets of $V$ of the form 
$a + Q^{k_{1}}L$. These are all disjoint and by the algorithm all of them 
must be chosen in the decomposition of $\mathcal{C}(V)$, and they all occur before
all the others. 
Thus $C_{1},\cdots,C_{r}$ and $C'_{1},\cdots,C'_{r}$ are $D_{1},\cdots,D_{r}$ 
in some order. Removing these and continuing in the same way the result is clear.
\hfill $\square$
\\

We have similar concepts in $\overline{L}$. For $W \subset \overline{L}$ we have 
the coset part $\mathcal{C}^{*}(W)$ of $W$ and $\mathcal{C}^{*}(W)$ can be written 
as a disjoint union of cosets in $W$. 
Let $\mathcal{C}^{*}(W) = \bigcup_{i} D_{i}$ where $D_{i},\, i=1,2,\cdots,$ 
are mutually disjoint cosets in $W$. 
We call $c^{*}(W) := \sum_i [\overline{L} : D_{i}]^{-1}$ 
the {\em total index} of $W$. This time we do not need to
be careful about the way in which the decomposition is obtained since the total
index is nothing else than the measure $\mu(\mathcal{C}^{*}(W))$ of 
$\mathcal{C}^{*}(W)$. 

Given an efficient decomposition $\mathcal{C}(V) = \bigcup_{i=1} C_{i}$ 
into disjont cosets in $L$, 
we define $\overline{\mathcal{C}}(V) := \bigcup_{i=1} \overline{C_{i}} 
\subset \overline{L}$. This is actually an open set in $\overline{L}$.
Since $[L:C] = [\overline{L}: \overline{C}]$
we see that $c(V) = c^*(\overline{\mathcal{C}}(V))$. In particular
it follows that the total index of any subset $V$ of $L$ is finite
and bounded by $\mu(\overline{\mathcal{C}}(V))$.

\begin{lemma}  \label{cInEqual}
For $X,Y \subset L$ and $X \subset Y$, and any decomposition
${\mathcal C}(X) = \bigcup_i C_{i}$ into disjoint cosets, 
$\sum_i c(C_i) \le c(Y)$. In particular, $c(X) \le c(Y)$.
\end{lemma}
{\sc proof}: Assume first that $Y$ is a single coset $C$.
Then 
\be
\sum_i c(C_i) = \sum_i c^*({\overline C}_i)
= \sum_i \mu ({\overline C}_i) \le \mu({\overline C}) = 
c^*({\overline C}) = c(C)\, ,
\ee
since the cosets remain distinct after closing them in ${\overline L}$.

In the general case, let 
$\mathcal{C}(Y) = {\stackrel{\cdot}{\bigcup}}_{j=1} C'_{j}$ 
be an efficient decomposition of $Y$. Since $X \subset Y$, 
each $C_{i} \subset Y$. In view of the remark above about
efficient decompositions, there is for each $i$ a unique $j$ for which 
$C_{i} \subset C'_{j}$. Thus we can arrange the $C_i$'s so that 
\be
\mathcal{C}(X) = \bigcup_{j=1}^\infty \bigcup_{i \in A_j} C_i
\ee
where $A_j:= \{i \,|\, C_i \subset C'_j \}$. Now
$\bigcup_{i \in A_j} C_i \subset C'_j$, so by the first part of the proof,
$\sum_{i \in A_j} c(C_i) \le c(C'_j)$. Finally 
\be
c(X) = \sum_j\sum_{i \in A_j}c(C_i) \le \sum c(C'_j) = c(Y) \, .
\ee
\hfill $\square$

\begin{lemma}
Let $U_{i},\, i = 1,\cdots,m,$ be disjoint point sets of the lattice 
$L$ in $\RR^{n}$. Let $\Lambda_{i} = \overline{U_{i}} \cap L$ 
and $\mathcal{C}(U_{i})$ be the coset part in $U_{i}$. 
Then 
$\bigcup_{i=1}^{m} (\Lambda_{i} \backslash U_{i}) \subset L 
\backslash \bigcup_{i=1}^{m} \mathcal{C}(U_{i})$, 
with equality if $L = \bigcup_{i=1}^{m} U_{i}$.
\end{lemma}
{\sc {proof}}
For $x \in \bigcup_{i=1}^{m} \mathcal{C}(U_{i})$ there is a coset 
$C \subset \mathcal{C}(U_{i})$ for which $x \in C \subset U_{i}$.
Let $C = a + Q^{k}L, a \in L$. Suppose $x$ is a limit point of $U_{j}$ 
in $\overline{L}$ for some $j \neq i$. Then, since $a + Q^{k}\overline{L}$ is an 
open neighborhood of $x$, $(a + Q^{k}\overline{L}) \cap U_{j} \neq \emptyset$
 ~ i.e. $(a + Q^{k}L) \cap U_{j} \neq \emptyset$.
But then $U_{i} \cap U_{j} \neq \emptyset$, contrary to the assumption. 
This means $x \not\in \bigcup_{i=1}^{m}(\Lambda_{i} \backslash U_{i})$, 
proving the first part.

Suppose that $L = \bigcup_{i=1}^{m}U_i$ and $x \in L$ but
$x \not\in \bigcup_{i=1}^{m} \mathcal{C}(U_{i})$. 
Then $x \in U_{i}$ for some $U_{i}$ but there is no coset in $U_{i}$ 
which contains $x$.
For any $k \in \ZZ_{+}, ~B_{k}(x) := x + Q^{k}\overline{L}$ is an open 
neighborhood of $x$ in $\overline{L}$ and $L \cap B_{k}(x) \not\subset U_{i}$,
 by assumption.
Since $L = \bigcup_{i=1}^{m} U_{i},
~ (L \cap B_{k}(x)) \cap U_{j} \neq \emptyset~$ for some $j \neq i$.
So we can choose $x_{k}^{j} \in (L \cap B_{k}(x)) \cap U_{j}$.
Then we get a sequence $\{x_{k}^{j}\}$ convergent to $x$ as $k \rightarrow \infty$.
Choosing a subsequence lying entirely in one $\Lambda_{j}$ shows that 
$x \in \Lambda_{j}$ for some $j \neq i$.
Since $x \in U_{i}$, and $U_{i},U_{j}$ are disjoint, 
$x \in \Lambda_{j} \backslash U_{j}$. \hfill $\square$

\begin{theorem}\label{mainCITheorem}
Let $U_{i},\, i = 1,\cdots,m,$ be disjoint nonempty point sets of the 
lattice $L$ in $\RR^{n}$. Let $\mathcal{C}(U_{i})$ be the coset part in 
$U_{i}$, $c(U_{i})$ the total index of $U_{i}$, and $W_i$ the closure
of $U_i$ in ${\overline L}$.
Then $\sum_{i=1}^{m} c(U_{i}) = 1$ if and only if the 
sets $U_{i},\, i = 1,\cdots,m,$ are regular weak model sets in the CPS(\rm{\ref{cps}}) 
and $\overline{L} = \bigcup_{i=1}^{m} W_{i}$.
\end{theorem}  
{\sc {proof}}\\
$(\Rightarrow)$
  Assume that $\sum_{i=1}^{m} c(U_{i}) = 1$. 
Let $U_{m+1} := L \backslash \bigcup_{i=1}^{m} U_{m}$. Using
Lemma \ref{cInEqual} and the fact that $c(L)=1$, we see that 
$c(U_{m+1}) = 0$ and  $\sum_{i=1}^{m+1} c(U_{i}) = 1$. For
this reason we can assume, in proving that the $U_i$
are weak model sets, that $\bigcup_{i=1}^{m}U_i = L$
in the first place.
 
For $j \ne k$ the cosets of $\mathcal{C}(U_{j})$ (of which
there may be none!) and those
of $\mathcal{C}(U_{k})$ are disjoint from one another, and the
same applies to $\overline{\mathcal{C}}(U_{j})$ and
$\overline{\mathcal{C}}(U_{k})$. Thus
\[
 \mu\left(\bigcup_{i=1}^{m} \overline{\mathcal{C}}(U_{i})\right) = 
\sum_{i=1}^{m}\mu( \overline{\mathcal{C}}(U_{i})) = 
\sum_{i=1}^{m} c(U_{i}) = 1, 
\]
and
\be \label{fullMeasure}
 \mu\left(\overline{L} \backslash 
(\bigcup_{i=1}^{m} \overline{\mathcal{C}}(U_{i}))\right)
 = 0.
\ee

Now note that 
$\partial W_j \cap \bigcup_{i=1}^{m}\overline{\mathcal{C}}(U_{i}) = \emptyset$
 for any $j$.
If not let $a \in \partial W_j \cap \overline{\mathcal{C}}(U_{k})$ for
some $k$. 
Since $\overline{\mathcal{C}}(U_{k}) \subset \, \stackrel{\circ}{W_{k}}$,
we see that $j \ne k$. But $a \in W_j$, so $a$ is a limit point of $U_j$
,and $\overline{\mathcal{C}}(U_{k})$ is an open neighborhood of $a$, 
so $U_j \cap \mathcal{C}(U_{k}) \ne \emptyset$. 
This violates the disjointness of the $U_i$'s.
We conclude that $\partial W_j \subset  \overline{L} \backslash 
(\bigcup_{i=1}^{m} \overline{\mathcal{C}}(U_{i}))$ and hence that
\be
\mu(\partial W_j) = 0\,,
\ee 
for all $j= 1,\dots,m.$
Note also
\begin{eqnarray}
{\Lambda_{i} \backslash U_{i}}
& \subseteq &{\bigcup_{j=1}^{m}(\Lambda_{j} \backslash U_{j})}\nonumber\\
& \subset &{L \backslash \bigcup_{j=1}^{m} \mathcal{C}(U_{j})}
 \quad \textrm{by Lemma 8} \nonumber \\ 
& = &{L \backslash \bigcup_{j=1}^{m} (\overline{\mathcal{C}}(U_{j}) \cap L)}
 = {L \backslash \bigcup_{j=1}^{m} \overline{\mathcal{C}}(U_{j})\,.}\nonumber
\end{eqnarray}
This shows that 
\be \label{measInequal}
\mu(\overline{\Lambda_{i} \backslash U_{i}}) \leq 
\mu\left(\overline{L \backslash (\bigcup_{i=1}^{m} 
\overline{\mathcal{C}}(U_{i}))}\right) \leq  
\mu\left(\overline{L} \backslash \bigcup_{i=1}^{m} 
\overline{\mathcal{C}}(U_{i})\right) = 0.
\ee
By Lemma \ref{3lemmas} (i) and (ii), 
$\stackrel{\circ}{W_{i}} \cap \stackrel{\circ}{W_{j}} = \emptyset$  
for all $ i \neq j$ and 
$~(\Lambda_{i} \backslash U_{i}) \subset \bigcup_{j=1}^{m} \partial W_{j}\,$ 
for all $\,i = 1,\cdots,m$. 

Using Lemma \ref{regModelSet}(i) we obtain that
the sets $U_{i}, i = 1,\cdots,m,$ are regular weak model sets in 
the CPS(\ref{cps}). 

Remark : Whenever $\stackrel{\circ}{W_i} \neq \emptyset$, 
$U_i$ is actually a regular model set.
 
Since $\bigcup_{i=1}^{m} \overline{\mathcal{C}}(U_{i}) 
\subset \bigcup_{i=1}^{m} W_{i},~ \mu(\bigcup_{i=1}^{m} W_{i}) = 1$.
Thus $\overline{L} \backslash \bigcup_{i=1}^{m} W_{i}$ is open of measure $0$ 
and $\overline{L} = \bigcup_{i=1}^{m} W_{i}$. This last argument
does not require that $\bigcup_{i=1}^{m}U_i = L$.
\\
\\
$(\Leftarrow)$
  Assume that $U_{i} = \Lambda(V_{i}) = V_{i} \cap L$ where 
$\overline{V_{i}} \backslash \stackrel{\circ}{V_{i}}$ has measure 
$0$ and  $\overline{L} = \bigcup_{i=1}^{m} W_{i}$. 
Thus $U_{i} \subset V_{i}$ and $W_{i} := \overline{U_{i}} \subset 
\overline{V_{i}}$. Since $L$ is dense in $\overline{L}$ and for 
$x \in \stackrel{\circ}{V_{i}}$ each ball around $x$ of radius $\epsilon > 0$ 
contains points of $\stackrel{\circ}{V_{i}} \cap L \subset U_{i}$, 
it follows that $\overline{U_{i}} \supset \stackrel{\circ}{V_{i}}$. 
This proves that $\stackrel{\circ}{V_{i}} \subset W_{i} \subset \overline{V_{i}}$. 
So $\mu(\stackrel{\circ}{V_{i}}) = \mu(W_{i})$ and 
$\mu(W_{i} \backslash \stackrel{\circ}{V_{i}}) = 0$.
Now 
\[\bigcup_{i=1}^{m} W_{i} = (\bigcup_{i=1}^{m} \stackrel{\circ}{V_{i}}) \cup
( \bigcup_{i=1}^{m}(W_{i} \backslash \stackrel{\circ}{V_{i}})).
\]
So $\mu(\bigcup_{i=1}^{m} W_{i}) = \mu(\bigcup_{i=1}^{m} \stackrel{\circ}{V_{i}})$.
 Also the disjointness of the $U_{i}$ gives
 $\stackrel{\circ}{V_{i}} \cap \stackrel{\circ}{V_{j}} = \emptyset$ 
for $i \neq j$ (since $L$ is dense in ${\overline L}$).
Finally
\[
1 = \mu(\bigcup_{i=1}^{m} W_{i}) = 
\mu(\bigcup_{i=1}^{m} \stackrel{\circ}{V_{i}}) = 
\sum_{i=1}^{m} \mu(\stackrel{\circ}{V_{i}}) = 
\sum_{i=1}^{m} c^{*}(\stackrel{\circ}{V_{i}}) \leq
\sum_{i=1}^{m} c(\stackrel{\circ}{V_{i}} \cap L) \leq \sum c(U_{i}) \leq 1\, .
\]  \hfill $\square$

\begin{coro}
 Let $(\tilde{U},\Phi)$ be a primitive substitution system with 
inflation $Q$ on the lattice $L$ in $\RR^{n}$.
  Suppose that PF-eigenvalue of the substitution matrix $S(\Phi)$ is equal to 
$|\det Q|$ and $L = \bigcup_{i=1}^{m} U_{i}$. 
Then $\sum_{i=1}^{m} c(U_{i}) = 1$, where $c(U_{i})$ is the total 
index of $U_{i}$ if and only if the sets $U_{i}, i = 1,\cdots,m,$ 
are model sets in CPS(\rm{\ref{cps}}).
\end{coro}  
{\sc {proof}}
Use Theorem \ref{mainPFtheorem} to determine that for all $i,
\stackrel{\circ}{W_i} \neq \emptyset$.
 Now use Theorem and Remark in the proof. 
\hfill $\square$

\section{Chair tiling}

\begin{figure}[ht]
\centerline{\epsfysize=30mm \epsfbox{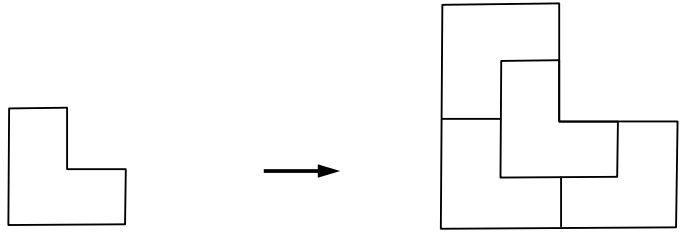}}
\caption{2-dimensional chair tiling inflation} \label{chair:1}
\end{figure}

The two dimensional chair tiling is generated by the inflation rule 
shown in Figure~\ref{chair:1}. There are $4$ orientations of the
chairs in any chair tiling. In \cite{padics} it was shown that the chair tiling
has an interpretation in terms of model sets based on the lattice $\ZZ^2$ and 
its 2-adic completion as internal space.

\begin{figure}[ht]
\centerline{\epsfysize=30mm \epsfbox{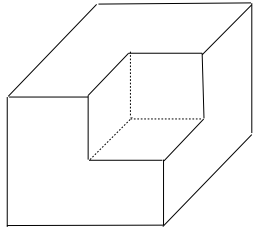}}
\caption{3-dimensional chair tile} \label{chair:3}
\end{figure}

In this section we generalize this result to the $n$-dimensional
chair tiling using the results of the last section
(see Figure \ref{chair:3} for an example of the $3$-dimensional chair). 
To make things clearer we begin with the case $n=2$. 

\subsection* { I. Chair tiling in $\RR^{2}$}

\begin{figure}[ht]
\centerline{\epsfysize=40mm \epsfbox{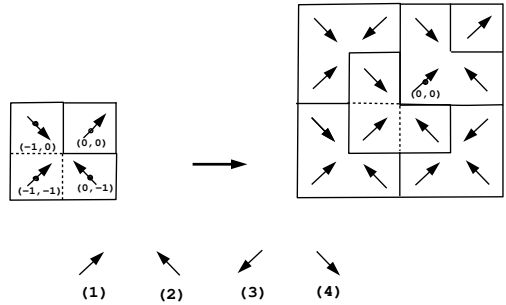}}
\caption{2-dimensional chair tiling substitution} \label{chair:2}
\end{figure}

The starting point is to replace each tile by $3$ oriented
squares. Figure~\ref{chair:2} shows the inflation rule, for one
chair, in terms of oriented squares.
The resulting tiling is a square tiling of the plane in which
each of the squares has one of $4$ orientations. The centre
points of each square form a square lattice which we identify
with $\ZZ^2$ by assigning coordinates as shown. 

Let $U_i$ be
the set of centre points corresponding to squares of orientation
$(i)$ as given in Figure~\ref{chair:2}.
We start out from a basic generating set 
$A_{2} := \{\,(x_{1},x_{2})\,|\,x_{i} \in \{0,-1\}\,\}$ and
determine the precise maps for the substitution rules of  
Figure~\ref{chair:2}.

Letting $e_{1} := (0,0),
 e_{2} := (1,0), e_{3} := (1,1), e_{4} := (0,1)\, ,$
these maps are defined as:
\[  
\begin{array}{ll}
f_{j,i} : U_{i} \rightarrow U_{j} ~~\mbox{by} ~(x,i) \mapsto (2x + e_{j},j) & \mbox{if}
~j \neq  i \pm 2\\
f_{i,i}^{(2)} : U_{i} \rightarrow U_{i} ~~\mbox{by} ~(x,i) \mapsto (2x + e_{j},i) &
\mbox{if} ~j = i \pm 2\,,\\
\end{array}
\]
where $i,j \in \{1,2,3,4\}, x \in \ZZ^{2},~i \pm 2 := \left\{ \begin{array}{ll}
                                            i+2 & ~\mbox{if}~ i \leq 2 \\
                                           i-2 & ~\mbox{if}~ i > 2\,. 
                                                    \end{array}
                                              \right.$\\
These are the maps of an affine substitution system $\Phi$.
In fact, if we define 
\[h_1:x \mapsto 2x+e_1, 
   ~ h_2:x \mapsto 2x + e_2,~ h_3:x \mapsto 2x +e_3,~ h_4:x \mapsto 2x +e_4,
\]
then
  \[ 
  \begin{array}{llll}
   f_{1,1} = h_1,& f_{1,2} = h_1,& f_{1,1}^{(2)} = h_3,& f_{1,4} = h_1\\
   f_{2,1} = h_2,& f_{2,2} = h_2,& f_{2,3} = h_2,& f_{2,2}^{(2)} = h_4\\
   f_{3,3}^{(2)} = h_1,& f_{3,2} = h_3,& f_{3,3} = h_3,& f_{3,4} = h_3\\ 
   f_{4,1} = h_4,& f_{4,4}^{(2)} = h_2,& f_{4,3} = h_4,& f_{4,4} = h_4\,,
 \end{array}   
  \]
and 
\[ \Phi = \left( \begin{array}{llll}
                  \{h_1,h_3\}& \{h_1\}& \{\} &\{h_1\}\\
                  \{h_2\} &\{h_2,h_4\}& \{h_2\}& \{\}\\
                  \{\}& \{h_3\}& \{h_3,h_1\}& \{h_3\}\\
                  \{h_4\}& \{\}& \{h_4\}& \{h_4,h_2\}\,
                  \end{array}
          \right)_{\,.}
\]

Inflating $A_{2}$ by the substitutions above we generate
the $4$
point sets $U_{i}, i = 1,2,3,4$. The precise description of $U_{i}$ is 
the following : 
\begin{eqnarray*}
{U_{1}}& = &
        {\bigcup_{k=0}^{\infty} \bigcup_{t=0}^{2^{k}-1} ((0,0)+2^{k} (2,0) +t(1,1) 
          +2^{k} \cdot 4 \ZZ^{2})}\\
       &&{ \cup~ \bigcup_{k=0}^{\infty} \bigcup_{t=0}^{2^k-1} ((0,0)+2^k (0,2) +t(1,1) 
          +2^k \cdot 4 \ZZ^2)}
     ~ \cup~ \bigcup_{t=-\infty}^{\infty} \{t(1,1)\} \\
{U_2}& = &
        {\bigcup_{k=0}^{\infty} \bigcup_{t=0}^{2^k-1} ((-1,0) +2^k (2,0) +t(-1,1) 
          +2^k \cdot 4 \ZZ^2)}\\
      && {\cup~ \bigcup_{k=0}^{\infty} \bigcup_{t=0}^{2^k-1} ((-1,0) +2^k (0,2)
          +t(-1,1) +2^k \cdot 4 \ZZ^2)}
     ~ \cup~ \bigcup_{t=0}^{\infty} \{(0,-1) +t(1,-1)\} \\
{U_3}& = &
        {\bigcup_{k=0}^{\infty} \bigcup_{t=0}^{2^k-1} ((-1,-1) +2^k (2,0) +t(-1,-1) 
          +2^k \cdot 4 \ZZ^2)} \\
        &&{\cup~ \bigcup_{k=0}^{\infty} \bigcup_{t=0}^{2^k-1} ((-1,-1) +2^k (0,2)
          +t(-1,-1) +2^k \cdot 4 \ZZ^2)} \\
{U_4}& = &
        {\bigcup_{k=0}^{\infty} \bigcup_{t=0}^{2^k-1} ((0,-1) +2^k (2,0) +t(1,-1) 
          +2^k \cdot 4 \ZZ^2)}\\
       && {\cup~ \bigcup_{k=0}^{\infty} \bigcup_{t=0}^{2^k-1} ((0,-1) +2^k (0,2) 
          +t(1,-1)+2^k \cdot 4 \ZZ^2)}
       ~ \cup~ \bigcup_{t=0}^{\infty} \{(-1,0) +t(-1,1)\}\,.
\end{eqnarray*}

Each of these decompositions is basically into cosets, with the exception
of three trailing sets in types $1, 2, 4$ which we will designate by
$V_1, V_2, V_4$ respectively.

We can prove the correctness of this as follows:

Let $U'_1 , U'_2 , U'_3, U'_4$ be the sets on the right hand sides respectively.
Note that 
\begin{enumerate}
  \item[(i)] The generating set $A_2$ is contained in $U'_i$ adequately, i.e. \\
              $~(0,0) \in U'_1, ~(0,-1) \in U'_2,~ (-1,-1) \in U'_1, 
               ~(-1,0) \in U'_4$.
  \item[(ii)] Claim that $U'_i \supset \bigcup_{j=1}^{4} \Phi_{ij} U'_j, i = 1,2,3,4$.

Check that for any $i$,
\begin{eqnarray*}
        { h_i(U'_i)}& \subset & { \bigcup_{\stackrel{j=1}{j \neq i,i \pm 2}}
            ^{4} \bigcup_{k=0}^{\infty} \bigcup_{t=0}^{2^k-1} \left( (-e_i)+
          2^{k+1}(2(e_i-e_j))+ 2t(e_{i \pm 2}-e_{i}) \right.}\\ 
                      & &{\left. ~~~~~~~~~~~~~~~~~~~~~~~~~+ 2^{k+1}
                        \cdot 4\ZZ^2 \right) ~{\cup ~V_i }}\\
                  & \subset & {U'_i}\\
         {h_{i \pm 2}(U'_i)}& \subset &{ \bigcup_{\stackrel{j=1}{j \neq i,i \pm 2}}^{4}
         \bigcup_{k=0}^{\infty} \bigcup_{t=0}^{2^k-1} \left( (-e_i)+2^{k+1}(2(e_i-e_j))+
                     (2t+1)(e_{i \pm 2}-e_{i}) \right.}\\
                               & &{\left.~~~~~~~~~~~~~~~~~~~~~~~~~+2^{k+1}
                                    \cdot 4\ZZ^2 \right) ~{\cup~V_i  } }\\
                                & \subset & {U'_i}\\
         {h_i(U'_l)}& \subset &{(-2e_l + e_i + 4\ZZ^2) \cup (-2e_{l \pm 2} + e_i+
           4\ZZ^2)}\\
                  & \subset & {U'_i\,, \mbox{where} ~l \neq i,i \pm 2,~l \in \{1,2,3,4\}}
  \end{eqnarray*}      

\item[(iii)] $~U'_i, i = 1,2,3,4,$ are all disjoint.\\
Indeed within each $U_i$, all the cosets and the  non-coset part
are clearly disjoint. 
And two cosets or non-coset sets chosen from $U'_i$ and $U'_j$, 
where $j \neq i,i \pm 2$, cannot intersect, since they are 
different by mod 2. Futhermore two of cosets or non-coset sets chosen from $U'_i$
and $U'_{i \pm 2}$ cannot intersect either, since for 
$a + 2^k\cdot4 \ZZ^2 \subset U'_i,~ b + 2^l\cdot4 \ZZ^2 \subset U'_{i \pm 2}$
with $ k \le l $, $~a - b \neq 0$  mod  $2^k \cdot 4 \ZZ^2$.
\end{enumerate}
Now since $U'_1,U'_2,U'_3,U'_4$ are generated from $A_2$ by $\Phi$, 
$~U_i \subset U'_i$  for all $i = 1,2,3,4$. Also from $\bigcup U_i = \ZZ^2$,
we get $\bigcup U'_i = \ZZ^2$ . Since all $U'_i, i = 1,2,3,4,$ are disjoint,
$U_i = U'_i$  for all $i = 1,2,3,4$.
   
Finally, for any $i = 1,2,3,4$,
\[ c(U_i) \geq 2 \cdot \sum_{k=0}^{\infty} \sum_{t=0}^{2^k-1} \frac{1}{(2^k\cdot4)^2}
 = 2 \cdot \sum_{k=0}^{\infty} \frac{2^k}{16 \cdot (2^k)^2} = \frac{1}{4}\,.
\]
Thus $\sum_{i=1}^{4} c(U_i) =1$.
Theorem \ref{mainCITheorem} shows that $U_i, i = 1,2,3,4,$ are regular model sets.

\subsection* {II. Chair tiling in $\RR^n$} 

In this section we are going to generalize the foregoing
to the $n$-dimensional chair tilings for all $n \ge 2$. The
$n$-chair is an $n$-cube with a corner taken out of it. The inflation
rule, which we spell out algebraically below, is geometrically 
the obvious generalization of the $2$-dimensional case.

We transform the geometry by replacing each chair by a $2^n -1$ oriented
cubes, as before, and coordinatize the lattice formed by the
centres of the cubes, starting from the basic generating set
$A_n := \{(x_1,\cdots,x_n) | x_i \in \{0,-1\}\}$. There are
$2^n$ orientations of cubes and hence $2^n$
{\em types} of points (but only $2^n -1$ of these types appear
in the starting set $A_n$).

For each $k \geq 0$ let $\beta(k)$ be the binary expansion
$\epsilon_o + \epsilon_1 2 + \epsilon_2 2^2 + \dots$ 
 of $k$,~$\epsilon_l \in \{0,1\}$. We define the basic orientation
vectors $e_1, \dots e_{2^n}$ by 
\[ 
e_i :=
\left\{ 
\begin{array}{ll}
(\epsilon_0, \dots , \epsilon_{n-1}) ~~ \textrm{ the binary digits of} ~ \beta(i-1)
~ & \textrm {if}~ i \le 2^{n-1} \, ,\\ 
(1, \dots 1) - e_{i - 2^{n-1}} & \textrm{if}~ i > 2^{n-1} \, .
\end{array}
\right. 
\]

We determine the sets $U_i$, 
$i = 1,\dots,2^n$, of all $i-$type points in $\ZZ^n$ from the 
points of the basic generating set $A_n$, using the inflation rules below.

The types of the points of $A_n$ are as follows:
for $x = (x_1,\dots,x_n) \in A_n$, 
\begin{itemize}
\item[] when $x_n = -1$ ;
  \begin{itemize}
   \item[] $~~x \in U_i$, for which 
           $\beta(i - 1) = (1,\dots,1) + x$,
  \end{itemize}
\item[] when $x_n = 0$ ;
   \begin{itemize}
    \item[] if $~ x = (0,\dots,0)$, $ ~~x \in U_1~~ $ 
    \item[] otherwise, $ ~~x \in U_{i + 2^{n-1}}$,
     for which $ ~\beta(i - 1) = (1,\dots,1) -((1,\dots,1) + x) ~~$.
   \end{itemize}
\end{itemize}
The idea of considering our vectors in the form $(1,\ldots,1) + x$ 
is to make it easy to compare them  with the  basic orientation vectors.

This conforms with what happens when $n = 2$: there
are $2^n -1$ types in the basic starting set that are in 
$2^{n-1} -1$ complementary pairs and $1$ pair of vectors
$(0, \dots, 0)$ and $(-1, \dots, -1)$ of the same type, namely
of type $1$.

Define 
\[  
\begin{array}{ll}
{f_{j,i} : U_{i} \rightarrow U_{j} ~~\mbox{by} ~(x,i) \mapsto (2x + e_{j},j)} &
{ \mbox{if}~j \neq i \pm 2^{n-1}}\\
{f_{i,i}^{(2)} : U_{i} \rightarrow U_{i} ~~\mbox{by} ~(x,i) \mapsto (2x + e_{j},i)}
 & {\mbox{if}~j = i \pm 2^{n-1}},
\end{array}
\]
where $i,j \in \{1,\cdots,2^n\},~~ x \in \ZZ^{n},~~i \pm 2^{n-1} := \left\{ \begin{array}{ll}
                    i+2^{n-1} & \mbox{if} ~i \leq 2^{n-1} \\
                    i-2^{n-1} & \mbox{if} ~i > 2^{n-1}\,. 
                   \end{array}
           \right.$\\
Let $\Phi$ be the matrix function system. Define $h_i : x \mapsto 2x+e_i,
     i \in \{1,\cdots,2^n\}$.
    \[ \Phi = \left( \begin{array}{c}
                      \{h_1,h_{1+2^{(n-1)}}\}~\{h_1\}~~~\cdots~\stackrel{\stackrel{1+2^{n-1}}
                       {\downarrow}}{\{\}}~ \cdots~~~\{h_1\}\\
                      \vdots\\
                      \{h_{2^n}\}~\{h_{2^n}\}~\cdots \stackrel{\stackrel{2^n-2^{n-1}}
                       {\downarrow}}{\{\}} \cdots~\{h_{2^n-2^{n-1}}, h_{2^n}\} 
                      \end{array}
              \right)_{\,.} \,
     \]\\
Inflating $A_n$ by the maps, we get the precise description of $U_i$ :
\[
 U_i = \bigcup_{\stackrel{j=1}{j \neq i,i \pm 2^{n-1} }}^{2^n} \bigcup_{k=0}^{\infty}
\bigcup_{t=0}^{2^k-1} \left( (-e_i)+2^k(2(e_i-e_j))+t(e_{i \pm 2^{n-1}}-e_{i})+
2^k \cdot 4\ZZ^n \right)~\bigcup~V_i\,,
\]
where 
\be
 V_i = \left\{ \begin{array}{ll}
                \bigcup_{t=-\infty}^{\infty} \{t(e_{1 \pm 2^{n-1}} - e_1)\} & \mbox{if}
                          ~i=1\\
                \bigcup_{t=0}^{\infty} \{t(e_i - e_{i \pm 2^{n-1}}) + 
                        (-e_{i \pm 2^{n-1}})\} & \mbox{if}~i \neq 1, 1 \pm 2^{n-1}  \\
                 \emptyset & \mbox{if} ~i=1 + 2^{n-1} \,.
                \end{array}
          \right.  \label{ch1}
\ee \\ 
The equalities can be proved in the same way as in the 2-dimensional case.\\
Let ${U'}_i$ be the set of the right hand side in (\ref{ch1}).
Note that
\begin{enumerate}
 \item[(i)]  The generating set $A_n$ is contained in $U'_i$ adequately, i.e.
            \begin{itemize}
               \item[] $e_1 \in U'_1$ \quad   if $ i = 1$
               \item[]  $-e_{i \pm 2^{n-1}} \in U'_i$ \quad if $~i \neq 1,
                1 \pm 2^{n-1}$
               \item[] $-e_{1 \pm 2^{n-1}} \in U'_1$ \quad if $~i = 1 \pm 2^{n-1}\,.$
            \end{itemize}
 \item[(ii)]  Claim that $U'_i \supset \bigcup_{j=1}^{2^n} \Phi_{ij} U'_j,
              i = 1,\cdots,2^n$.\\
 Indeed for $i \in \{1,\cdots,2^n\}$
 \begin{eqnarray*}
        { h_i(U'_i)}& \subset & {\bigcup_{\stackrel{j=1}{j \neq i,i \pm 2^{n-1}}}
            ^{2^n} \bigcup_{k=0}^{\infty} \bigcup_{t=0}^{2^k-1} \left( (-e_i)+
          2^{k+1}(2(e_i-e_j))+ 2t(e_{i \pm 2^{n-1}}-e_{i}) \right.}\\
                        & &{\left. ~~~~~~~~~~~~~~~~~~~~~~~~~+ 2^{k+1} 
                            \cdot 4\ZZ^n \right) ~{\cup ~V_i } }\\
                  & \subset & {U'_i}\\
         {h_{i \pm 2^{n-1}}(U'_i)}& \subset &{ \bigcup_{\stackrel{j=1}
             {j \neq i, i \pm 2^{n-1}}}^{2^n}
         \bigcup_{k=0}^{\infty} \bigcup_{t=0}^{2^k-1} \left( (-e_i)+2^{k+1}(2(e_i-e_j))+
          (2t+1)(e_{i \pm 2^{n-1}}-e_{i}) \right.}\\
                                & &{\left.~~~~~~~~~~~~~~~~~~~~~~~~~ + 2^{k+1} 
                                    \cdot 4\ZZ^n \right) {\cup~V_i  } }\\
                                 & \subset &{U'_i}\\
         {h_i(U'_l)}& \subset &{(-2e_l + e_i + 4\ZZ^n) \cup (-2e_{l \pm 2^{n-1}} + e_i+
           4\ZZ^n)}\\
                  & \subset & {U'_i\,,~\mbox{where} ~l \neq i,~i \pm 2^{n-1} ,~ l \in \{1,\cdots,2^n\}\,.}
  \end{eqnarray*}        
 \item[(iii)] $U'_i, i = 1,\cdots,2^n,$ are disjoint.\\
Indeed all cosets and a non-coset set in each $U'_i$ are all disjoint.
And two cosets or non-coset sets chosen from $U'_i$ and $U'_j$, 
where $~j \neq i,i \pm 2^{n-1}$, cannot intersect, since they are 
different by mod 2. Futhermore two of cosets or non-coset sets chosen from $U'_i$
and $U'_{i \pm 2^{n-1}}$ cannot intersect either, since for 
$a + 2^k\cdot4 \ZZ^{n} \subset U'_i,~ b + 2^l\cdot4 \ZZ^n \subset U'_{i \pm 2^{n-1}}$
with $ k \leq l $,  $~a - b \neq 0$  mod  $2^k \cdot 4 \ZZ^n$.
\end{enumerate}
Now since $U'_i,i = 1,\cdots,2^n,$ are generated from $A_n$ by $\Phi$, 
$~U_i \subset U'_i$ for all $i = 1,\cdots,2^n$. Also from $\bigcup_{i=1}^{2^n} U_i = \ZZ^n$,
$\bigcup_{i=1}^{2^n} U'_i = \ZZ^n$. Since all $U'_i, i = 1,\cdots,2^n,$ are disjoint,
$U_i = U'_i$  for all $i = 1,\cdots,2^n$.\\
   
For any $i = 1,\cdots,2^n$,
\[ c(U_i) \geq (2^n-2) \cdot \sum_{k=0}^{\infty} \sum_{t=0}^{2^k-1} \frac{1}{(2^k\cdot4)^n}
 = (2^n-2) \cdot \sum_{k=0}^{\infty} \frac{2^k}{2^{2n} \cdot (2^k)^n} = 
\frac{1}{2^n}.
\]
Thus $\sum_{i=1}^{2^n} c(U_i) =1$.
Theorem \ref{mainCITheorem} shows that $U_i, i = 1,\cdots,2^n,$ are regular model sets.

To get a model set interpretation of the chair tiling itself we proceed as 
follows.
We observe that every arrow points to the inner corner of exactly one chair. 
Let us label each chair by its inner corner point which is at the tip of exactly
$2^n-1$ arrows. These corner points give us $2^n$ sets $X_1,\ldots,X_{2^n}$ 
according to the type, and all lie in the shift $L' = (\frac{1}{2},\ldots,
\frac{1}{2}) + \ZZ^n$ of our lattice $\ZZ^n$.
Let $f_i, i=1,\ldots,2^n$, be $(\frac{1}{2},\ldots,\frac{1}{2}) - e_i$ respectively. Then $U_i + f_i$ is 
the set of tips of all arrows of type $i$ and $U_i + f_i = L' \cap (V_i + f_i)\,,$ 
for some $V_i \subset \ZZ_2^n$ for which $\overline{V_i}$ compact, $\stackrel{\circ}{V_i} \ne \emptyset$
and $\mu(\partial V_i) = 0$.
Now \[X_i = L' \cap (\bigcap_{j \neq i \pm 2^{n-1}}(V_j + f_j))\] which is the 
required regular
model set description of $X_i$, since 
\[\partial(\bigcap_{j \neq i \pm 2^{n-1}}
(V_j + f_j)) \subset \bigcup_{j \neq i \pm 2^{n-1}}\partial(V_j + f_j)
\] and 
$\mu(\partial(V_j + f_j)) = 0$ for all $j = 1,\ldots,2^n$.\\

From this result we can show that if we mark each chair with a single
point in a consistent way, then the set of points obtained from all the
chairs of any one type also forms a regular model set, and hence a pure point
diffractive set.

\section*{Acknowledgment} It is a pleasure to thank Martin Schlottmann
for his encouragement and in particular for his help in Theorem \ref{mainPFtheorem}. 
We are also indebted to Michael Baake for his interest and insights into this work.

\end{document}